\documentclass{amsart}

\usepackage{amsmath,amsthm,amssymb,amsfonts}
\usepackage{bbm}
\usepackage[mathscr]{eucal}
\usepackage[english]{babel}
\usepackage{graphicx}

\allowdisplaybreaks[1]

\newtheorem{thm}{Theorem}[section]

\newtheorem{lem}[thm]{Lemma}
\newtheorem{prop}[thm]{Proposition}
\newtheorem{rem}[thm]{Remark}

\newtheorem{defn}[thm]{Definition}

\theoremstyle{remark}

\numberwithin{equation}{section}

\def\BB{{\mathbb B}}
\def\CC{{\mathbb C}}
\def\NN{{\mathbb N}}
\def\RR{{\mathbb R}}
\def\SS{{\mathbb S}}

\def\cB{\mathcal{B}}

\def\cD{\mathcal{D}}
\def\cE{\mathcal{E}}
\def\cF{\mathcal{F}}

\def\cH{\mathcal{H}}

\def\cM{\mathcal{M}}

\def\cP{\mathcal{P}}

\def\cS{\mathcal{S}}
\def\cV{\mathcal{V}}
\def\cX{\mathcal{X}}

\def\fD{{\mathfrak D}}

\def\sD{D}

\def\eps{{\varepsilon}}
\def\supp{\operatorname{supp}}

\def\QQ{{\cQ}}

\def\cc{{\mathfrak c}}

\def\fD{{\mathfrak D}}
\def\KK{{\mathscr{K}}}

\def\PP{{\mathsf {P}}}

\def\ZZZ{Z}
\def\tZZZ{{\overline{Z}}}
\def\YYY{Y}
\def\tYYY{{\overline{Y}}}
\def\omeg{\omega}
\def\tP{{\widetilde{P}}}

\def\LLL{{\mathscr{L}}}
\def\LLV{{\mathscr{L}_V}}

\def\QQ{{Q}}
\newcounter{rea}
\begin{document}

\title[Gaussian bounds for the heat kernel associated to the PSWFs]
{Gaussian bounds for the heat kernel associated\\
to prolate spheroidal wave functions\\
with applications}

\author[A. Bonami]{Aline Bonami}
\address{Institut Denis Poisson, (UMR 7013), Universit\'e d’Orl\'eans, Universit\'e
de Tours \& CNRS,  45067 Orl\'eans, France}
\email{aline.bonami@univ-orleans.fr}

\author[G. Kerkyacharian]{Gerard Kerkyacharian}
\address{University Paris Diderot-Paris 7, LPMA, Paris, France}
\email{kerk@math.univ-paris-diderot.fr}

\author[P. Petrushev]{Pencho Petrushev}
\address{Department of Mathematics\\
University of South Carolina\\
Columbia, USA}
\email{pencho@math.sc.edu}

\thanks{The second author has been supported ANR Grant Basics.}
\thanks{The third author has been supported by NSF Grant DMS-1714369.}
%\thanks{Corresponding author: Pencho Petrushev, E-mail: pencho@math.sc.edu}

\begin{abstract}
Gaussian upper and lower bounds and H\"older continuity are established for the heat kernel
associated to the prolate spheroidal wave functions (PSWFs) of order zero.
These results are obtained by application of a general perturbation principle
using the fact that the PSWF operator is a perturbation of the Legendre operator.
Consequently, the Gaussian bounds and H\"older inequality for the PSWF heat kernel follow from the ones in the Legendre case.
As an application of the general perturbation principle, we also establish Gaussian bounds
for the heat kernels associated to generalized univariate PSWFs and PSWFs on the unit ball in $\RR^d$.
Further, we develop the related to the PSWFs of order zero smooth functional calculus,
which in turn is the necessary groundwork in developing the theory
of Besov and Triebel-Lizorkin spaces associated to the PSWFs.
One of our main results on Besov and Triebel-Lizorkin spaces associated to the PSWFs
asserts that they are the same as the Besov and Triebel-Lizorkin spaces generated by the Legendre operator.
\end{abstract}

%\date{March 8, 2020}
%\date{August 18, 2020}
%\date{December 26, 2020}
%\date{January 6, 2021}
%\date{August 25, 2021}

\maketitle

\section{Introduction}\label{sec:introduction}

The prolate spheroidal wave functions (PSWFs) of order zero are widely used in physics (see e.g. \cite{F, MF}),
signal processing, and numerical analysis (see e.g. \cite{HL,ORX}).
Their importance is rooted in the fact that the PSWFs are band-limited functions
that are eigenfunctions of two important operators:
\begin{equation}\label{def:Lc}
L_\cc f(x):= -\frac{d}{dx}\Big[(1-x^2)\frac{df}{dx}(x)\Big] + \cc^2x^2f(x),
\quad x\in (-1, 1),\; \cc>0,
\end{equation}
and
\begin{equation}\label{def-Fc}
F_\cc f(x):= \int_{-1}^1f(t)e^{i\cc xt}dt
\end{equation}
that commute.
As is customary, we denote by $\psi_0, \psi_1, \dots$
the eigenfunctions of the operators $L_\cc$ and $F_\cc$ (the PSWFs),
and by $0<\chi_0<\chi_1 <\cdots$ the eigenvalues of $L_\cc$
and by $\lambda_0, \lambda_1, \dots$,
where $|\lambda_0| > |\lambda_1| > \cdots >0$,
and $\lambda_n=i^n|\lambda_n|$,
the eigenvalues of $F_\cc$.
Thus,
\begin{equation}\label{eigenf-L}
L_\cc\psi_n=\chi_n\psi_n,
\quad n=0, 1, \dots,
\end{equation}
\begin{equation}\label{eigenf-F}
F_\cc\psi_n=\lambda_n\psi_n,
\quad n=0, 1, \dots.
\end{equation}
Of course, the eigenfunctions $\{\psi_n\}$ and eigenvalues $\{\chi_n\}$ and $\{\lambda_n\}$
depend on the constant $\cc$.
We assume that $\cc > 0$ is arbitrary but fixed throughout.

It is easy to see that the PSWFs $\{\psi_n\}$ are also eigenfunctions of the operator
\begin{equation}\label{def:Qc}
Q_\cc f(x):= \frac{2\pi}{\cc}\big(F_\cc^*F_\cc\big)f(x)
=\frac{1}{\pi}\int_{-1}^1 \frac{\sin \cc(x-y)}{x-y}f(y) dy,
\end{equation}
that is,
\begin{equation}\label{eigenf-Q}
Q_\cc\psi_n=\mu_n\psi_n,
\quad \mu_n = \frac{\cc}{2\pi} |\lambda_n|^2,
\quad n=0, 1, \dots.
\end{equation}
Observe that the operator $L_\cc$ is essentially self-adjoint and positive.

We assume that the PSWFs $\{\psi_n\}$ are normalized in $L^2[-1, 1]$ and as is well known
\begin{equation}\label{orthog}
\langle \psi_n, \psi_k \rangle := \int_{-1}^1 \psi_n(x)\psi_k(x)dx = \delta_{nk},
\quad \overline{\psi_n(x)}= \psi_n(x).
\end{equation}
It is also well known that the PSWFs $\{\psi_n\}_{n\ge 0}$ form and orthonormal basis for $L^2[-1, 1]$.

It should be pointed out that the PSWFs $\{\psi_n\}$ extend analytically as functions on $\RR$ that are band-limited
with band limit $\cc$, that is,
$\supp \widehat{\psi}_n \subset [-\cc, \cc]$.
Also, they are orthogonal on~$\RR$:
$\int_\RR \psi_n(x)\psi_k(x)dx=0$ if $n\ne k$.
Furthermore, the PSWFs $\{\psi_n\}$ form an orthogonal basis for
the space of band-limited functions on $\RR$ with band limit $\cc$.

We refer the reader to \cite{HL, ORX} for details and references on the PSWFs.

One of our goals in this article is to establish Gaussian upper and lower bounds for the heat kernel
associated to the operator $L_\cc$ as well as its Lipschitz continuity.
The parameter $\cc$ of the operator $L_\cc$ can be given value $\cc=0$.
Then the functions $\{\psi_n\}$ are the normalized Legendre polynomials.
Gaussian upper and lower bounds for the associated heat kernel are established in \cite{CKP}.
The Gaussian bounds for the heat kernel associated to the operator $L_\cc$ when $\cc>0$
are derived from the respective estimates for the heat kernel associated to $L_0$
by an application of a general principle that relates semigroups
associated to a self-adjoint operator and its perturbation and their kernels.

As an application of the Gaussian bounds for the heat kernel associated to the operator $L_\cc$
%these results
we develop a smooth functional calculus,
which in turn enables us to develop the theory of Besov and Triebel-Lizorkin spaces on $[-1, 1]$
based on the PSWFs $\{\psi_n\}$ in the spirit of Frazier and Jawerth \cite{FJ1, FJ2, FJW}.
In developing this theory we use the approach and ideas from \cite{CKP, KP}.

We will use the remainder of this introduction to give
a more detailed outline of the contents of this article.

\subsection{The PSWF heat kernel}\label{subsec:PSWF-HK}

The operator $L_\cc$ is (essentially) self-adjoint and positive
and hence it generates a semigroup $\exp(-tL_\cc)$, $t>0$.
It is easy to see that $\exp(-tL_\cc)$ is an integral operator with (heat) kernel
$p_t(x, y)$ defined by
\begin{equation}\label{def-heat-k}
p_t(x, y):=\sum_{n=0}^\infty e^{-t\chi_n}\psi_n(x)\psi_n(y),
\quad x, y\in [-1, 1].
\end{equation}
As will be explained it is natural to consider the interval $[-1, 1]$ equipped with
the usual Lebesgue measure $d\mu(x) = dx$
and the distance $\rho(x, y):= |\arccos x - \arccos y|$.
We use the standard notation $B(x, r):= \{y\in [-1, 1]:\rho(x, y)<r\}$
for the ``balls'' on $[-1,1]$ and
set $V(x, r):=\mu(B(x, r))$.

Gaussian lower and upper bounds for the PSWF heat kernel $p_t(x, y)$ are established (Theorem~\ref{thm:main}):
%There exist constants $c_1, c_2, c_3, c_4 >0$  such that
For any $x,y\in [-1, 1]$ and $t>0$
\begin{equation}\label{Gauss-0}
\frac{c_1\exp\big(- \frac{\rho(x,y)^2}{c_2t}\big)}{\big[V(x, \sqrt t)V(y, \sqrt t)\big]^{1/2}}
\le  p_t(x, y)
\le \frac{c_3\exp\big(- \frac{\rho(x,y)^2}{c_4t}\big)}{\big[V(x, \sqrt t)V(y, \sqrt t)\big]^{1/2}}.
\end{equation}
The following Lipschitz continuity of $p_t(x, y)$ is also obtained (Theorem~\ref{thm:Holder-PSWF}):
\begin{equation}\label{Holder-0}
|p_t(x, y)-p_t(x', y)|
\le c\frac{\rho(x, x')}{\sqrt{t}}
\cdot\frac{\exp\big(- \frac{\rho(x,y)^2}{ct}\big)}{\big[V(x, \sqrt t)V(y, \sqrt t)\big]^{1/2}},
\end{equation}
if $\rho(x, x')\le \sqrt{t}$.

To establish these estimates %the Gaussian bounds \eqref{Gauss-0} and Lipschitz continuity \eqref{Holder-0}
we use a general perturbation principle that we outline next.

\subsection{Perturbation of self-adjoint operators and associated semigrous}\label{subsec:principle}

A~general perturbation principle is our tool
in establishing the Gaussian bounds and H\"{o}lder continuity of the heat kernels associated
to PSWFs, generalized PSWFs and PSWFs on the ball.
Given a measure space, we assume that $\YYY$ and $\ZZZ$ are two self-adjoint positive operator,
and $\YYY$ is a perturbation of $\ZZZ$ of the form:
$$
\YYY f := \ZZZ f + Af
\quad\hbox{with}\quad
Af(x) := V(x)f(x),
$$
where $V\ge 0$ is a bounded function. %$L^2$-function.
Denote by $S_t$ and $T_t$ the one parameter semigroups generated by $\ZZZ$ and $\YYY$, respectively.
Also, let $\cE(f, g):=\langle \ZZZ f, g\rangle$
be the associated to $\ZZZ$ quadratic form.
Assuming that the extension $\overline{\cE}$ of $\cE$ is a Dirichlet form (\S\ref{sec:abstract})
and that the operator $S_t: L^1\to L^\infty$ is bounded with kernel $\KK_{S_t}(x, y)$
we show that
$T_t$, $t>0$, is an integral operator with kernel
$\KK_{T_t} \in L^\infty$
satisfying
\begin{equation}\label{KS-KT}
e^{-t\|V\|_\infty}\KK_{S_t}(x, y)\le \KK_{T_t}(x, y) \le \KK_{S_t}(x, y)
\quad\hbox{almost everywhere}, \;\; \forall t>0.
\end{equation}
These inequalities show that $\KK_{T_t}(x, y)$ inherits the bounds for $\KK_{S_t}(x, y)$.
In particular, $\KK_{T_t}$ inherits the Gaussian bounds of $\KK_{S_t}$ if such bounds exist,
which is the case in our context.
Inequalities similar to the ones from \eqref{KS-KT}, under adequate assumptions,
have been developed by different authors,
see e.g. \cite[Theorem  2.21]{Ouhabaz} and \cite[Theorem 3.6]{ArDe}.

To establish inequalities \eqref{KS-KT} we utilize a well known variation of parameter formula (see \eqref{Tt-St})
to represent $T_t$ in terms of $S_t$, as given in textbooks,
see e.g. Corollary 1.7 in \cite{Engel-Nagel}
or Theorem~11.4.1 and identity (11.12), p. 341, in \cite{Davies-2}.
This formula also allows to show that
the kernel $\KK_{T_t} $ inherits the H\"older continuity of  $\KK_{S_t}$ (see \eqref{Holder-0})
in the case of a doubling metric measure space whenever Gaussian bounds for $\KK_{S_t}(x, y)$ exist.
For details see Section~\ref{sec:abstract}.

\subsection{Functional calculus induced by the PSWFs of order zero}

Our further development
heavily depends on the almost exponential space localization of
the kernels
$\KK_{F(\delta\sqrt{L_\cc})}(x, y)$ of operators $F(\delta\sqrt{L_\cc})$, $\delta>0$, for smooth functions $F$.
These are kernels of the form
\begin{equation}\label{ker-fLc}
\KK_{F(\delta\sqrt{L_\cc})}(x, y)= \sum_{n=0}^\infty F(\delta\sqrt{\chi_n})\psi_n(x)\psi_n(y).
%\quad x, y\in [-1, 1]. %, \; \delta>0.
\end{equation}
As is shown in Theorem~\ref{thm:ker-local-S}
if the function $F$ is even, real-valued, and belongs to the Schwartz class $\cS(\RR)$ on $\RR$,
then for any $\sigma>0$ there exists a constant $c_\sigma>0$ such that
\begin{equation*}
|\KK_{F(\delta\sqrt{L_\cc})}(x,y)|
\le \frac{c_\sigma\Big(1+\frac{\rho(x, y)}{\delta}\Big)^{-\sigma}}{\big[V(x, \delta)V(y, \delta)\big]^{1/2}},
\quad x, y\in [-1, 1],\; \delta >0.
\end{equation*}
Furthermore, the kernel $\KK_{F(\delta\sqrt{L_\cc})}(x,y)$ is Lipschitz continuous (see Theorem~\ref{thm:ker-lip}).
The proof of these estimates relies on the finite speed propagation property of the relate to $L_\cc$ wave equation
(see Proposition~\ref{prop:f-speed}),
a consequence of the Gaussian localization of the PSWF heat kernel.
%%
%For kernels as above we also prove the following important localization estimates %(Lemma~\ref{lem:L0-Lc}):
%(Theorem~\ref{thm:L0-Lc}):
%For any $\sigma>0$ and $k\in \NN$
%\begin{equation*}
%|L_{0, x}^k \KK_{F(\delta\sqrt{L_\cc})}(x, y)|, |L_{\cc, x}^k \KK_{F(\delta\sqrt{L_0})}(x, y)|
%\le \frac{c(\sigma, k)\delta^{-2k}}{V(x, \delta)\Big(1+\frac{\rho(x, y)}{\delta}\Big)^{\sigma}},
%\quad 0<\delta \le 1,
%\end{equation*}
%where the operators $L_0^k$ and $L_\cc^k$ are acting on the kernels as functions of $x$.

\subsection{PSWF Besov and Triebel-Lizorkin spaces} %\label{subsec:B-F-spaces}

One of the aims of this article is to develop the Besov and Triebel-Lizorkin spaces associated to the PSWFs.
In doing so we have been guided by ideas from \cite{CKP, KP}.
The Besov and Triebel-Lizorkin spaces are in general spaces of distributions.
We introduce distributions associated with the operators $L_0$ and $L_\cc$ in Subsection~\ref{subsec:distributions}.
It is easy to show that the distributions in both cases are the same.

Using kernels as in \eqref{ker-fLc} we introduce two kinds of Besov spaces associated to the operators $L_0$, $L_\cc$:
(a) Classical Besov spaces $B^s_{pq}(L_0)$, $\cB^s_{pq}(L_\cc)$, and
(b) Non-classical Besov spaces $\widetilde{B}^s_{pq}(L_0)$, $\widetilde{\cB}^s_{pq}(L_\cc)$
with full sets of parameters: $s\in\RR$ and $0<p, q\le \infty$.

One of our main results asserts that
$B^s_{pq}(L_0)=\cB^s_{pq}(L_\cc)$
and
$\widetilde{B}^s_{pq}(L_0)=\widetilde{\cB}^s_{pq}(L_\cc)$
with equivalent norms (see Section~\ref{sec:B-F-spaces}).

We also use kernels as in \eqref{ker-fLc} to introduce
classical and non-classical Triebel-Lizorkin spaces
$F^s_{pq}(L_0)$, $\cF^s_{pq}(L_\cc)$ and
$\widetilde{F}^s_{pq}(L_0)$, $\widetilde{\cF}^s_{pq}(L_\cc)$,
associated to the operators $L_0$ and $L_\cc$,
$s\in\RR$, $0<p<\infty$, and $0<q\le \infty$.
Furthermore, we show that
$F^s_{pq}(L_0)=\cF^s_{pq}(L_\cc)$
and
$\widetilde{F}^s_{pq}(L_0)=\widetilde{\cF}^s_{pq}(L_\cc)$
with equivalent norms.

As is well known the PSWFs $\{\psi_n\}$ are closely related to the Legendre polynomials.
The above results show that in terms of Besov and Triebel-Lizorkin spaces on $[-1, 1]$
the PSWFs generate the same spaces as the Legendre polynomials.
%
%The construction of frames based on the PSWFs $\{\psi_n\}$ will be discussed
%in Section~\ref{sec:fremes}.

\subsection{Generalized PSWFs and PSWFs on the ball}

We consider two additional important cases where our general result, discussed in \S\ref{subsec:principle},
produce Gaussian upper and lower bounds for heat kernels.

The Jacobi operator $\LLL$ is defined by
\begin{equation*}
\LLL f(x):=-\frac{1}{\omega(x)}\frac{d}{dx}\Big[\omega(x)a(x)\frac{df}{dx}(x)\Big],
\quad D(\LLL) := C^2[-1, 1],
\end{equation*}
where
\begin{equation*}
\omega(x)=\omega_{\alpha, \beta}(x):= (1-x)^\alpha(1+x)^\beta,
\quad
\alpha, \beta >-1,
\quad a(x):=(1-x^2).
\end{equation*}
Consider now a {\bf perturbation of the operator} $\LLL$ of the form
\begin{equation*}
\LLV f(x):= \LLL f(x)+ V(x)f(x),
\quad \quad D(\LLV) := C^2[-1, 1],
\end{equation*}
where $V\in C[-1, 1]$ and $V\ge 0$.
In this setting it is natural to consider the distance
$\rho(x, y):=|\arccos x - \arccos y|$, just as in the case of the PSWFs.

As an application of the general principle described in \S\ref{subsec:principle} we establish
Gaussian upper and lower bounds (similar to \eqref{Gauss-0}) for the heat kernel associated
to the operator $\LLV$, see Section~\ref{sec:Jacobi}.
In the particular case when $\alpha = \beta$ and $V(x)=\cc^2 x^2$
this leads to Gaussian bounds for the heat kernel associated to the generalized PSWFs introduced in \cite{ZLWZ}.
The H\"{o}lder continuity of these kernels is also established.

\smallskip

The prolate spheroidal wave functions on the unit ball $\BB^d$ in $\RR^d$
were introduced by D. Slepian in \cite{Slep-1}.
They are defined as the eigenfunctions of the following two commuting operators:
\begin{equation*}%\label{def-L-ball-intro}
L_\cc f(x):=-\sum_{i=1}^{d}\partial_i^2f(x) + \sum_{i=1}^d \sum_{j=1}^d x_ix_j  \partial_{i}\partial_jf(x)
+ (n+2\gamma)\sum_{i=1}^{d} x_i\partial_if(x) + c^2\|x\|^2f(x)
\end{equation*}
and
\begin{equation*}%\label{def-Fc-B}
F_\cc f(x) :=\int_{\BB^d} f(u)e^{-i\cc x\cdot u}\omega_\gamma(u)du,
\end{equation*}
where $\omega_\gamma(x):= (1-\|x\|^2)^{\gamma-1/2}$, $\gamma >-1/2$, is a weight function
with $\|x\|$ being the Euclidean norm of $x\in \RR^d$
and $\cc>0$ is a fixed constant.
The natural distance here is defined by
$\rho(x,y) := \arccos \big(x\cdot y + \sqrt{1-\| x\|^2}\sqrt{1-\| y\|^2}\big)$.

Yet another application of the general principle, described in \S\ref{subsec:principle}, allows to
establish Gaussian upper and lower bounds (as in \eqref{Gauss-0}) for the heat kernel
associated to the operator $L_\cc$, see Section~\ref{sec:PSWF-ball}.

\subsection*{Notation}

We will denote by $\langle f, g \rangle$ the inner product in the respective Hilbert space
and use the abbreviated notation $\|\cdot\|_p:=\|\cdot \|_{L^p}$.
The norm of a bounded operator $L: L^p\to L^q$ will be denoted by $\|L\|_{p\to q}$.
As usual $\NN$ will denote the set of all positive integers
and $\NN_0:=\NN\cup\{0\}$.
The notation $a\vee b:= \max\{a, b\}$ and $a\wedge b:= \min\{a, b\}$ will also be used.
Positive constants will be denoted by $c, c_1, c_2, c', \dots$
and they may vary at every occurrence;
$a\sim b$ stands for $c_1\le a/b\le c_2$.

\section{Perturbation of self-adjoint operators and associated semigroups}\label{sec:abstract}

As already alluded to in the introduction the proof of the Gaussian bounds for the PSWF heat kernel (see \eqref{Gauss-0})
as well as the heat kernels associated to the generalized PSWFs and the PSWFs on the ball
will depend on an abstract result relating semigroups associated to
a self-adjoint operator and its perturbation and their kernels.
We next present this result.
In developing our general theory we will use some basic facts from the theory of
self-adjoint operators and the related quadratic forms and one-parameter semigroups
that can be found in e.g. \cite{A, Davies-1, Davies-2, Fukushima}.

We will consider two self-adjoint operators $\ZZZ$ and $\YYY$ in a Hilbert space
with the second being a perturbation of the first,
and the associated to them one-parameter semigroups.

\subsection{The operator $\ZZZ$ and the associated semigroup}\label{subsec:Oper-Z}

We assume that $(X, \mu)$ is a $\sigma$-finite measure space and
consider the real Hilbert space $\cH$ consisting of all real-valued functions in $L^2(X, \mu)$.

Let $\ZZZ$ be a positive self-adjoint operator on $\cH=L^2(X, \mu)$
with domain $D(\ZZZ)$ that is dense in $\cH$, i.e. $\overline{D(\ZZZ)}=\cH$.
This means that $\ZZZ$ is a positive symmetric operator whose closure
$\tZZZ$ (with domain $D(\tZZZ)$) is positive and self-adjoint.
From Spectral Theory it follows that the operator $\tZZZ$ generates
a self-adjoint one parameter contraction semigroup $S_t= e^{-t \tZZZ}$, $t>0$, on $\cH$:
\begin{equation*}
\|S_tf\|_2 \le \|f\|_2, \quad \forall f\in \cH. %f\in L^2(X, \mu).
\end{equation*}
Consider the associated to $\ZZZ$ symmetric non-negative quadratic form $\cE$:
\begin{equation*}
\cE(f, g):=\langle \ZZZ f, g\rangle,
\qquad
\cE(f, f)=\langle \ZZZ f, f\rangle \ge 0,
\end{equation*}
with domain $D(\cE)=D(L)$.
As is well known $\cE$ is closable
and denoting by $\overline{\cE}$ and $D(\overline{\cE})$ its closure and domain
one has
\begin{equation}\label{quad-form}
\overline{\cE}(f, g) = \big\langle \tZZZ^{\frac{1}{2}}f, \tZZZ^{\frac{1}{2}}g\big\rangle
\quad\hbox{and}\quad
D(\overline{\cE}) = D\big(\tZZZ^{\frac{1}{2}}\big).
\end{equation}

Our next assumption is that $\overline{\cE}$ is a Dirichlet form on $\cH$.
Namely, we stipulate the following condition:
\begin{equation}\label{BD-cond}
f\in D(\overline{\cE})
\quad \Longrightarrow\quad
(f\vee 0)\wedge 1 \in D(\overline{\cE}),
\quad
\overline{\cE}\big((f\vee 0)\wedge 1,(f\vee 0)\wedge 1\big) \le \overline{\cE}(f, f).
\end{equation}

The point is that because $\overline{\cE}$ is a Dirichlet form on $\cH$,
then the semigroup $S_t$ is submarkovian:
If $0\le f\le 1$ and $f\in L^2(X, \mu)$, then
$0\le S_tf\le 1$.
%This readily implies that $S_tf\ge 0$ if $f\in L^2(X, \mu)$ and $f\ge 0$.

In addition, we assume that
\begin{equation}\label{St-bound}
\|S_tf\|_\infty \le \Upsilon(t)\|f\|_1,
\quad \forall f\in L^1(X, \mu),
\end{equation}
for some function $\Upsilon(t)>0$.
Then as is well known (see \cite{DS}, Theorem 6, p. 503)
$S_t$ is an integral operator with kernel
$\KK_{S_t} \in L^\infty$, i.e.
\begin{equation*}
S_tf(x)= \int_X \KK_{S_t}(x, y)f(y)d\mu(y), \quad \forall f\in L^1(X, \mu),
\end{equation*}
and
$\|S_t\|_{1\to\infty} = \|\KK_{S_t}\|_\infty \le \Upsilon(t)$.
The fact that the semigroup $S_t$ is submarkovian implies that
$S_tf\ge 0$ if $f\ge 0$ and $f\in L^1(X, \mu)$,
which leads to the conclusion that
\begin{equation}\label{positive-K}
\KK_{S_t}(x, y) \ge 0
\quad\hbox{for a.a.} \;\; x, y \in X\times X.
\end{equation}

\subsection{Perturbation $\YYY$ of the operator $\ZZZ$ and the associated semigroup}\label{subsec:perturbation}

We now consider a second operator $\YYY$ defined by
\begin{equation}\label{def-Y}
\YYY f := \ZZZ f + Af,\quad f\in D(\YYY) := D(\ZZZ),
\end{equation}
where the operator $A$ is defined by
\begin{equation}\label{def-A}
Af(x) := V(x)f(x)
\quad\hbox{with}\quad
%V\in L^\infty(X, \mu) \cap L^2(X, \mu), \quad V\ge 0.
V\in L^\infty(X, \mu), \quad V\ge 0.
\end{equation}
Clearly,
\begin{align*}
\langle \YYY f, g\rangle = \langle f, \YYY g\rangle,
\;\; \forall f, g\in D(\YYY)
\quad\hbox{and}\quad
\langle \YYY f, f\rangle \ge 0, \;\; \forall f\in D(\YYY).
\end{align*}
Therefore,
$\YYY$ is a positive self-adjoint operator on $L^2(X, \mu)$
and its closure is given by
$\tYYY=\tZZZ+A$ with domain $D(\tYYY):=D(\tZZZ)$.
Let $T_t=e^{-t\tYYY}$ be the associated to $\tYYY$ one parameter contraction semigroup on $L^2(X, \mu)$.
The associated to $\YYY$ symmetric non-negative quadratic form $\cE_\YYY$
is defined by
\begin{equation*}
\cE_\YYY(f, g):= \langle \YYY f, g \rangle = \cE(f, g) + \int_X V(x)f(x)g(x) d\mu(x),
\quad f, g\in D(\cE)=D(\ZZZ).
\end{equation*}
For its closure $\overline{\cE_\YYY}$ we have
\begin{equation*}
\overline{\cE_\YYY}(f, g)= \overline{\cE}(f, g) + \int_X V(x)f(x)g(x) d\mu(x),
\quad f, g\in D(\overline{\cE_\YYY}):= D(\overline{\cE}).
\end{equation*}
We claim that $\overline{\cE_\YYY}$ is a Dirichlet form on $L^2(X, \mu)$.
Indeed, by our assumption $\overline{\cE}$ is a Dirichlet form on $L^2(X, \mu)$,
i.e. \eqref{BD-cond} is valid.
Using this we have for any $f\in D(\overline{\cE_\YYY}):= D(\overline{\cE})$,
with $h:=(f\vee 0)\wedge 1$,
\begin{align*}
\overline{\cE_\YYY}(h, h)
&= \overline{\cE}(h, h) + \int_X V(x)h^2(x) d\mu(x)
\\
&\le \overline{\cE}(f, f)+ \int_X V(x)f^2(x) d\mu(x)
= \overline{\cE_\YYY}(f, f).
\end{align*}
Hence, $\overline{\cE_\YYY}$ is a Dirichlet form.
Just as before this implies that the semigroup $T_t$ is submarkovian:
If $0\le f\le 1$ and $f\in L^2(X, \mu)$, then
$0\le T_tf\le 1$.
Consequently, if $f\ge 0$ and $f\in L^2(X, \mu)$, then
$T_tf\ge 0$.

The following theorem will play a key role in establishing Gaussian bounds for heat kernels.

\begin{thm}\label{thm:abstract}
In the setting described above,
$T_t$, $t>0$, is an integral operator with kernel
$\KK_{T_t} \in L^\infty$, i.e.
\begin{equation*}
T_tf(x)= \int_X \KK_{T_t}(x, y)f(y)d\mu(y), \quad \forall f\in L^1(X, \mu),
\end{equation*}
satisfying
\begin{equation}\label{est-K-K}
e^{-t\|V\|_\infty}\KK_{S_t}(x, y)\le \KK_{T_t}(x, y) \le \KK_{S_t}(x, y)
\quad\hbox{for a.a.}\; x, y \in X, \; \forall t>0.
\end{equation}
\end{thm}

\begin{proof}
Using the terminology from \cite{Davies-2} and \cite{Engel-Nagel},
the semigroup $S_t$ is generated by the operator $-\tZZZ$ and $T_t$ is generated by $-\tYYY$.
Then $-\tYYY=-\tZZZ-A$ leads to the following identity:
\begin{equation}\label{Tt-St}
T_tf = S_tf-\int_0^t S_{t-s}AT_s f ds, \quad t>0, \; f\in L^1(X, \mu),
\end{equation}
see Theorem~11.4.1 and identity (11.12), p. 341, in \cite{Davies-2}
or Corollary 1.7 in \cite{Engel-Nagel}.

Assume $f\in L^1(X, \mu)$ and $f\ge 0$.
Since $Af(x):= V(x)f(x)$, $V\ge 0$, and $S_t$, $T_t$ are positivity preserving we have
$AT_sf \ge 0$ and $S_{t-s}AT_s f \ge 0$.
This coupled with \eqref{Tt-St} leads to
\begin{equation}\label{Tt-le-St}
0\le T_tf = S_tf-\int_0^t S_{t-s}AT_s f ds \le S_tf
\end{equation}
and using \eqref{St-bound} it follows that
\begin{equation}\label{Tt-le-St-norm}
\|T_tf\|_\infty \le \|S_tf\|_\infty \le \Upsilon(t)\|f\|_1.
\end{equation}

Let $f\in L^1(X, \mu)$.
Clearly, $f=f_+-f_-$, where $f_+:= f\vee 0$ and $f_-:= -(f\wedge 0)$;
obviously $f_+\ge 0$ and $f_-\ge 0$.
Then using the above we get
\begin{align*}
\|T_tf\|_\infty &\le \|T_tf_+\|_\infty + \|T_tf_-\|_\infty
\le \|S_tf_+\|_\infty + \|S_tf_-\|_\infty
\\
&\le \Upsilon(t)\|f_+\|_1 +\Upsilon(t)\|f_-\|_1
= \Upsilon(t)\|f\|_1.
\end{align*}
As before (\cite{DS}, Theorem 6, p. 503) we conclude that
$T_t$ is an integral operator with kernel
$\KK_{T_t} \in L^\infty$
such that
$\|T_t\|_{1\to\infty} = \|\KK_{T_t}\|_\infty \le \Upsilon(t)$.
From this and \eqref{Tt-le-St} it follows that
if $f\in L^1(X, \mu)$ and $f\ge 0$, then for almost all $x\in X$
\begin{align*}
0 \le T_tf(x)= \int_X \KK_{T_t}(x, y)f(y)d\mu(y)
\le S_tf(x)= \int_X \KK_{S_t}(x, y)f(y)d\mu(y),
\end{align*}
which readily implies
\begin{equation*}
\KK_{T_t}(x, y) \le \KK_{S_t}(x, y)
\quad\hbox{for a.a.}\; x, y\in X\times X.
\end{equation*}
The proof of the upper bound estimate in \eqref{est-K-K} is complete.

For the estimate in the other direction we use \eqref{Tt-St}, \eqref{Tt-le-St},
and the fact that $S_t$, $T_t$ preserve positivity
to obtain for any $f\in L^1(X, \mu)$, $f\ge 0$,
\begin{align*}
S_tf &= T_tf+\int_0^t S_{t-s}AT_s f ds
\le T_tf+ \|V\|_\infty\int_0^t S_{t-s}T_s f ds
\\
&\le T_tf+ \|V\|_\infty\int_0^t S_{t-s}S_s f ds
= T_tf + t\|V\|_\infty S_tf.
\end{align*}
Hence,
\begin{align*}
(1-t\|V\|_\infty) S_tf \le T_tf
\quad\Longrightarrow\quad
\Big(1-\frac{t\|V\|_\infty}{n}\Big) S_{t/n}f \le T_{t/n}f,
\quad t>0,\; n\in\NN.
\end{align*}
Using the fact that $S_t$, $T_t$ are semigroups (e.g. $S_{t_1+t_2}=S_{t_1}S_{t_2}$) this yields
\begin{equation*}
\Big(1-\frac{t\|V\|_\infty}{n}\Big)^n S_tf \le T_tf
\end{equation*}
and letting $n\to \infty$ we arrive at
\begin{equation*}
e^{-t\|V\|_\infty}S_tf \le T_t f,
\end{equation*}
which implies the lower bound estimate in \eqref{est-K-K}.
\end{proof}

We conclude this subsection by some bibliographic comments.
The semigroup $S_t$ dominates the semigroup $T_t$ following the usual terminology, see for instance \cite{Ouhabaz}.
The fact that such a domination translates into an inequality on kernels is given in \cite{Davies-2}, Problem 11.4.3.
It is an easy consequence of the Lebesgue differentiation theorem.

\subsection{Gaussian bounds and H\"{o}lder continuity of the kernels \boldmath $\KK_{S_t}$ and $\KK_{T_t}$}\label{subsec:Holder}

We assume that the quadratic form $\overline{\cE}$ from \eqref{quad-form} is
{\em regular and strictly local}.
Then $(\overline{\cE}, D(\overline{\cE}))$
is called a regular strictly local Dirichlet space.
In addition, we assume that there is a distance $\rho(x,y)$ on $X$ that is equivalent
to the intrinsic distance generated by the quadratic form $\overline{\cE}$.
For the definitions and details, see e.g. \cite{CKP}, \cite{GS}, and the references therein.
As in the introduction, we define
\begin{equation*}
V(x, r):= \mu(B(x, r)),
\quad
B(x, r):= \{y\in X: \rho(x, y)<r\}.
\end{equation*}
We assume that the heat kernel $\KK_{S_t}(x, y)$ has the following Gaussian bounds:
\begin{equation}\label{Gauss-SZ}
\frac{c_1\exp\big(- \frac{\rho(x,y)^2}{c_2t}\big)}{\big[V(x, \sqrt t)V(y, \sqrt t)\big]^{1/2}}
\le \KK_{S_t}(x, y)
\le \frac{c_3\exp\big(- \frac{\rho(x,y)^2}{c_4t}\big)}{\big[V(x, \sqrt t)V(y, \sqrt t)\big]^{1/2}},
%\quad \forall x,y\in X, t>0,
\end{equation}
for all $x,y\in X$, $t>0$.
The following upper and lower bounds for $\KK_{T_t}$ follow from above and Theorem~\ref{thm:abstract}:
%From above and Theorem~\ref{thm:abstract} it follows that estimates are valid for $\KK_{T_t}$:
\begin{equation}\label{Gauss-TZ}
\exp(-t\|V\|_\infty)\frac{c_1\exp\big(- \frac{\rho(x,y)^2}{c_2t}\big)}{\big[V(x, \sqrt t)V(y, \sqrt t)\big]^{1/2}}
\le \KK_{T_t}(x, y)
\le \frac{c_3\exp\big(- \frac{\rho(x,y)^2}{c_4t}\big)}{\big[V(x, \sqrt t)V(y, \sqrt t)\big]^{1/2}}.
\end{equation}

As is well known, the above Gaussian bounds of the heat kernel $\KK_{S_t}(x, y)$
are equivalent to the fact that:
(a) The measure $\mu$ has the doubling property and
(b) The respective local scale-invariant Poincar\'{e} inequality is valid.
See e.g. \cite[Theorem 2.31]{GS}.
Recall that the measure $\mu$ has the doubling property if
\begin{equation}\label{doubling-X}
0< V(x, 2r) \le cV(x, r), \quad \forall x\in X, r>0,
\end{equation}
where $c>0$ is a constant.

It is also well known that, in the present setting, the Gaussian bounds \eqref{Gauss-SZ}
of the kernel $\KK_{S_t}(x, y)$ imply its H\"{o}lder continuity:
There exists a constant $0<\alpha\le 1$ such that
\begin{equation}\label{S-Holder}
|\KK_{S_t}(x, y)-\KK_{S_t}(x', y)|
\le c_5\Big(\frac{\rho(x, x')}{\sqrt{t}}\Big)^\alpha
\frac{\exp\big(- \frac{\rho(x,y)^2}{c_6t}\big)}{\big[V(x, \sqrt t)V(y, \sqrt t)\big]^{1/2}},
\end{equation}
for all $x, x', y \in X$ and $t>0$, whenever $\rho(x, x')\le \sqrt{t}$.
We refer the reader to \cite[Theorem 2.32]{GS} and also to \cite{BCF, HS, SC, DP}
for details on the H\"{o}lder continuity of heat kernels.

We next show that $\KK_{T_t}(x, y)$
inherits the H\"{o}lder continuity \eqref{S-Holder} of $\KK_{S_t}(x, y)$.

\begin{thm}\label{thm:Holder}
In the setting of this section,
for any $x, x', y \in X$ and $t>0$
\begin{equation}\label{T-Holder}
|\KK_{T_t}(x, y)-\KK_{T_t}(x', y)|
\le c_7(1+t\|V\|_\infty) \Big(\frac{\rho(x, x')}{\sqrt{t}}\Big)^\alpha
\frac{\exp\big(- \frac{\rho(x,y)^2}{c_8t}\big)}{\big[V(x, \sqrt t)V(y, \sqrt t)\big]^{1/2}},
\end{equation}
whenever $\rho(x, x')\le \sqrt{t}$.
%where $c_7, c_8>0$ are constants.
\end{thm}

\begin{proof}
From identity \eqref{Tt-St} it follows that
\begin{equation*}
\KK_{T_t}(x, y) = \KK_{S_t}(x, y) -\int_0^t \int_X \KK_{S_{t-s}}(x, z)V(z)\KK_{T_t}(z, y) dz ds,
\end{equation*}
implying
\begin{multline}\label{KT-KT}
|\KK_{T_t}(x, y)-\KK_{T_t}(x', y)| \le |\KK_{S_t}(x, y) - \KK_{S_t}(x', y)|
\\
+\int_0^t \int_X |\KK_{S_{t-s}}(x, z)-\KK_{S_{t-s}}(x', z)|V(z)\KK_{T_t}(z, y) dz ds
\\
\le |\KK_{S_t}(x, y) - \KK_{S_t}(x', y)|
+\int_0^t \int_X |\KK_{S_{t-s}}(x, z)-\KK_{S_{t-s}}(x', z)|V(z)\KK_{S_t}(z, y) dz ds,
\end{multline}
where for the last inequality we used \eqref{est-K-K}.

The doubling property \eqref{doubling-X} readily implies that there exist constants $d, c_0>0$ such that
\begin{equation}\label{doubling-X2}
V(x, \lambda r) \le c_0\lambda^d V(y, r), \quad \forall x\in X, r>0, \lambda \ge 1.
\end{equation}
In turn, just as in \eqref{double-3} this yields
\begin{equation}\label{V-V}
V(x, r) \le c_0\Big(1+\frac{\rho(x, y)}{r}\Big)^d V(y, r), \quad \forall x, y\in X, r>0.
\end{equation}

Estimates \eqref{Gauss-SZ}, \eqref{S-Holder} together with \eqref{V-V} yield
\begin{equation}\label{Gauss-SZ2}
0\le \KK_{S_t}(x, y)
\le \frac{c_9\exp\big(- \frac{\rho(x,y)^2}{c_{10}t}\big)}{V(y, \sqrt t)},
\quad \forall x,y\in X, t>0,
\end{equation}
and
\begin{equation}\label{S-Holder2}
|\KK_{S_t}(x, y)-\KK_{S_t}(x', y)|
\le c_{11}\Big(\frac{\rho(x, x')}{\sqrt{t}}\Big)^\alpha
\frac{\exp\big(- \frac{\rho(x,y)^2}{c_{10}t}\big)}{V(x, \sqrt t)},
\end{equation}
for all $x, x', y \in X$ and $t>0$, whenever $\rho(x, x')\le \sqrt{t}$.
Clearly, \eqref{S-Holder2} holds with the roles of $x$ and $x'$ interchanged.
This along with \eqref{Gauss-SZ2} and \eqref{S-Holder2} leads to
\begin{equation*}%\label{S-Holder22}
|\KK_{S_t}(x, y)-\KK_{S_t}(x', y)|
\le c\Big(\frac{\rho(x, x')}{\sqrt{t}}\Big)^\alpha
\left(\frac{\exp\big(- \frac{\rho(x,y)^2}{c_{10}t}\big)}{V(x, \sqrt t)}
+\frac{\exp\big(- \frac{\rho(x',y)^2}{c_{10}t}\big)}{V(x', \sqrt t)}
\right),
\end{equation*}
for all $x, x', y \in X$ and $t>0$ (without the restriction $\rho(x, x')\le \sqrt{t}$).
In turn, this estimate along with \eqref{KT-KT}, \eqref{S-Holder}, and \eqref{Gauss-SZ2} imply that
whenever $\rho(x, x')\le \sqrt{t}$
\begin{align}\label{KTKT}
&|\KK_{T_t}(x, y)-\KK_{T_t}(x', y)|
\le c_5\Big(\frac{\rho(x, x')}{\sqrt{t}}\Big)^\alpha
\frac{\exp\big(- \frac{\rho(x,y)^2}{c_6t}\big)}{[V(x, \sqrt t)V(y, \sqrt t)]^{1/2}}
\\
& + c \|V\|_\infty\rho(x, x')^\alpha
\int_0^t (t-s)^{-\alpha/2}\int_X
\frac{\exp\big(- \frac{\rho(x,z)^2}{c_{10}(t-s)}\big)}{V(x, \sqrt{t-s})}
\cdot\frac{\exp\big(- \frac{\rho(z,y)^2}{c_{10}s}\big)}{V(y, \sqrt s)} dzds \nonumber
\\
& + c \|V\|_\infty\rho(x, x')^\alpha
\int_0^t (t-s)^{-\alpha/2}\int_X
\frac{\exp\big(- \frac{\rho(x',z)^2}{c_{10}(t-s)}\big)}{V(x', \sqrt{t-s})}
\cdot\frac{\exp\big(- \frac{\rho(z,y)^2}{c_{10}s}\big)}{V(y, \sqrt s)} dzds. \nonumber
\end{align}

We next estimate the integrals above over $X$.

\begin{lem}\label{lem:int-X}
For any  $x,y\in X$, $t>0$, and $0<s<t$ we have
\begin{equation}\label{est-int-X}
\int_X
\frac{\exp\big(- \frac{\rho(x,z)^2}{c(t-s)}\big)}{V(x, \sqrt{t-s})}
\cdot\frac{\exp\big(- \frac{\rho(z,y)^2}{cs}\big)}{V(y, \sqrt s)} dz
\le  \frac{c''\exp\big(- \frac{\rho(x,y)^2}{c't}\big)}{[V(x, \sqrt t)V(y, \sqrt t)]^{1/2}},
\end{equation}
where $c, c', c''>0$ are constants.
\end{lem}
\begin{proof}
We consider two cases.
First, let $0<s\le t/2$.
Considering two subcases $\rho(x,z)\le \rho(x,y)/2$ and $\rho(x,z)> \rho(x,y)/2$ it follows that
\begin{equation*}
\frac{\rho(x,z)^2}{c(t-s)}+ \frac{\rho(z,y)^2}{cs}
\ge \frac{\rho(x,y)^2}{8ct}+ \frac{\rho(z,y)^2}{2cs}.
\end{equation*}
Using this and that
$V(x, \sqrt{t-s}) \ge V(x, \sqrt{t/2}) \ge c V(x, \sqrt{t})$, in light of \eqref{doubling-X2},
we obtain
\begin{align*}
&\int_X
\frac{\exp\big(- \frac{\rho(x,z)^2}{c(t-s)}\big)}{V(x, \sqrt{t-s})}
\cdot\frac{\exp\big(- \frac{\rho(z,y)^2}{cs}\big)}{V(y, \sqrt s)} dz
\le c\frac{\exp\big(- \frac{\rho(x,y)^2}{8ct}\big)}{V(x, \sqrt{t})}
\int_X \frac{\exp\big(- \frac{\rho(z,y)^2}{2cs}\big)}{V(y, \sqrt s)} dz
\\
& \le c'\frac{\exp\big(- \frac{\rho(x,y)^2}{8ct}\big)}{V(x, \sqrt{t})V(y, \sqrt s)}
\int_X \Big(1+\frac{\rho(z,y)}{\sqrt{s}}\Big)^{-d-1} dz
\le c''\frac{\exp\big(- \frac{\rho(x,y)^2}{8ct}\big)}{V(x, \sqrt{t})}.
\end{align*}
Clearly, \eqref{est-int-X} follows from this and \eqref{V-V}.
Above we used the following simple inequality (see \cite[Lemma 2.3]{CKP}):
\begin{equation}\label{basic-ineq}
\int_X \Big(1+\frac{\rho(y,z)}{\delta}\Big)^{-d-1}dz \le cV(y, \delta),
\quad \forall y\in X, \delta >0,
\end{equation}
where the constant $d>0$ is from \eqref{doubling-X2}.
The proof of \eqref{est-int-X} in the case when $t/2<s\le t$ is similar (symmetric) and will be omitted.
\end{proof}

Now, we use \eqref{KTKT}, Lemma~\ref{lem:int-X}, and
that $\int_0^t(t-s)^{-\alpha/2}ds =\frac{t^{1-\alpha/2}}{1-\alpha/2}$
to obtain
\begin{align}\label{KT-KT-F}
&|\KK_{T_t}(x, y)-\KK_{T_t}(x', y)|
\\
&\le c(1+t\|V\|_\infty)\Big(\frac{\rho(x, x')}{\sqrt{t}}\Big)^\alpha
\Big(\frac{\exp\big(- \frac{\rho(x,y)^2}{ct}\big)}{[V(x, \sqrt t)V(y, \sqrt t)]^{1/2}}
+\frac{\exp\big(- \frac{\rho(x',y)^2}{ct}\big)}{[V(x', \sqrt t)V(y, \sqrt t)]^{1/2}}\Big), \nonumber
\end{align}
whenever $\rho(x, x')\le \sqrt{t}$.
Note that because $\rho(x, x')\le \sqrt{t}$ it follows from \eqref{V-V} that
$V(x,\sqrt{t}) \le c_02^d V(x',\sqrt{t})$
and obviously
$$
\frac{\rho(x,y)^2}{t} \le \Big(\frac{\rho(x, x')+\rho(x',y)}{\sqrt{t}}\Big)^2
\le 2\Big(1+\frac{\rho(x',y)^2}{t}\Big).
$$
We use these inequalities in \eqref{KT-KT-F} to conclude that estimate \eqref{T-Holder} is valid
whenever $\rho(x, x')\le \sqrt{t}$.
\end{proof}

\section{Gaussian bounds for the heat kernel associated to PSWFs}\label{sec:PSWF-0}

In this section we present some additional properties of the PSWFs $\{\psi_n\}$
of order zero, introduced in the introduction,
and establish Gaussian upper and lower bounds and the Lipschitz continuity for the associated heat kernel.

\subsection{The geometry of the space. Legendre polynomials}\label{subsec:prelim}

The fact that the prolate spheroidal wave functions are closely related to the Legandre polynomials
leads to the conclusion that it is natural to consider %in our setting
the interval $[-1,1]$ equipped with the Lebesgue measure $d\mu(x)=dx$
and the distance
\begin{equation}\label{dist}
\rho(x, y):= |\arccos x - \arccos y|.
\end{equation}
This is most clearly reflected in the Gaussian upper and lower bounds \eqref{Gauss-L0}
for the heat kernel generated by the Legendre polynomials.

Recall our notation
$B(x, r):= \{y\in [-1, 1]:\rho(x, y)<r\}$
for the ``balls" on $[-1,1]$.
As is well known (see e.g. \cite{CKP}) in the natural range $0<r\le \pi$
\begin{equation}\label{mes-B}
V(x, r):=\mu(B(x, r)) \sim r(1-x+r^2)^{1/2}(1+x+r^2)^{1/2}
\sim r(\sqrt{1-x^2}+r)
\end{equation}
and
$V(x,r) = V(x, \pi)=2$ if $r>\pi$.
This readily implies that the measure $d\mu=dx$ has the doubling property:
\begin{equation}\label{doubling}
V(x, 2r) \le c V(x, r), \quad x\in [-1, 1], \; r>0.
\end{equation}
We will also use that
\begin{equation}\label{mes-B-inv}
V(x, r)^{-1}\sim\min\left\{ \frac 1{r^2}, \frac 1{r\sqrt{1-x^2}}\right\}, \quad 0<r\le \pi.
\end{equation}
Observe that from \eqref{mes-B} it follows that
\begin{equation}\label{double-2}
V(x, \lambda r)\le c_0\lambda^2V(x, r),
\quad x\in [-1, 1],\; r>0,\; \lambda \ge 1,
\end{equation}
where $c_0>0$ is an absolute constant.
Since $B(x, r)\subset B(y, \rho(x, y)+r)$, \eqref{double-2} yields the following useful inequality
\begin{equation}\label{double-3}
V(x, r) \le c_0\Big(1+\frac{\rho(x, y)}{r}\Big)^2V(y, r),
\quad x, y\in [-1, 1],\; r>0.
\end{equation}

To establish Gaussian bounds for the PSWF heat kernel $p_t(x, y)$, defined in \eqref{def-heat-k},
and its Lipschitz continuity
we will use the fact that the operator $L_\cc$
is a perturbation of the Legendre operator
\begin{equation}\label{def:L0}
L_0 f(x):= -\frac{d}{dx}\Big[(1-x^2)\frac{df}{dx}(x)\Big],
\quad x\in (-1, 1).
\end{equation}
As is well known the Legendre polynomials $\{\bar{P}_n\}$ are eigenfunctions of the operator $L_0$.
More precisely,
\begin{equation}\label{P-eigen-f}
L_0\bar{P}_n = n(n+1) \bar{P}_n, \quad n\ge 0.
\end{equation}
We will assume that the Legendre polynomials $\{\bar{P}_n\}$ are normalized in $L^2[-1, 1]$.
The kernel $\KK_{\exp (-tL_0)}(x,y)$ %$\exp (-tL_0)(x,y)$
of the semigroup $\exp (-tL_0)$ takes the form
\begin{equation}\label{heat-ker-L0}
\KK_{\exp (-tL_0)}(x,y) =\sum_{n\ge 0} e^{-tn(n+1)} \bar{P}_n(x)\bar{P}_n(y).
\end{equation}

As is shown in \cite[Theorem~7.2]{CKP} (see also \cite[Theorem~5.1]{KPX1})
the (heat) kernel $\KK_{\exp (-tL_0)}(x,y)$ has Gaussian upper and lower bounds:
%
%There exist constants $c_1, c_2, c_3, c_4 >0$  such that for all $x,y\in [-1, 1]$ and $t>0$
\begin{equation}\label{Gauss-L0}
\frac{c_1\exp\big(- \frac{\rho(x,y)^2}{c_2t}\big)}{\big[V(x, \sqrt t)V(y, \sqrt t)\big]^{1/2}}
\le  \KK_{\exp (-tL_0)}(x,y)
\le \frac{c_3\exp\big(- \frac{\rho(x,y)^2}{c_4t}\big)}{\big[V(x, \sqrt t)V(y, \sqrt t)\big]^{1/2}}.
\end{equation}
In fact, the above Gaussian bounds are established in \cite{CKP}, \cite{KPX1}
in the more general case of the heat kernel associated to the Jacobi operator.

\subsection{Properties of the PSWFs}\label{subsec:prop-PSWF}

As we already alluded to in the introduction
we consider the parameter $\cc>0$ in the definition of the PSWFs selected arbitrarily but fixed.
We will need the following estimates on the eigenvalues $\{\chi_n\}$ (see e.g. inequalities (3.4) in \cite{ORX}):
\begin{equation}\label{eigen-v}
n(n+1)<\chi_n< n(n+1)+\cc^2.
\end{equation}
The PSWFs $\{\psi_n\}$ can be regarded as a perturbation of the Legendre polynomials on $[-1, 1]$.
%Denote by $\bar{P}_n$ the $n$th degree Legendre polynomial normalized in $L^2[-1,1]$,
%that is, $\|\bar{P}_n\|_{L^2}=1$.

As is well known the Legendre polynomials $\{\bar{P}_n\}_{n\ge 0}$ form an orthonormal basis for $L^2[-1,1]$.
Note also that
\begin{align}
\|\bar{P}_n\|_{L^\infty[-1, 1]} &= \bar{P}_n(1)=\sqrt{n+1/2}, \quad\hbox{and} \label{LP-est-1}
\\
|\bar{P}_n(x)|&\le \frac{(3/\pi)^{1/2}}{(1-x^2)^{1/4}},
\quad x\in (-1, 1). \label{LP-est-2}
\end{align}
The first equalities above are well known \cite{Sz};
the second inequality follows from inequality (7.3.8) in \cite{Sz} and \eqref{LP-est-1}.

We next compare $\psi_n$ with $\bar{P}_n$.

\begin{prop}\label{prop:PSWF}
We have, for all $n\ge 0$,
\begin{equation}\label{eigen-f}
\|\psi_n-\bar{P}_n\|_{L^\infty[-1, 1]}\le \frac{2\cc^2}{\sqrt{n+1/2}},
\end{equation}
and
\begin{equation}\label{eigen-f-2}
|\psi_n(x)-\bar{P}_n(x)|\le \frac{c}{(n+1)(1-x^2)^{1/4}},
\quad x\in (-1, 1),
\end{equation}
where the constant $c>0$ depends only on $\cc$.
Consequently,
\begin{equation}\label{norm-psi}
 \|\psi_n(x)\|_{L^\infty[-1, 1]} \le c(n+1)^{1/2}
\end{equation}
and
\begin{equation}\label{est-psi}
|\psi_n(x)|\le \frac{c}{(1-x^2)^{1/4}},
\quad x\in (-1, 1).
\end{equation}
\end{prop}
\begin{proof}
Estimate \eqref{eigen-f} is a part of \cite[Proposition~5]{BK}.
Estimate \eqref{eigen-f-2} can be deduced from the proof of \cite[Proposition~5]{BK}, p. 43, as follows.
In this proof, it is shown that
there exists a constant $A =A(n)$ such that for $x\in [-1, 1]$
\begin{equation}\label{est-PSWF}
|\psi_n(x)-A\bar{P}_n(x)|\le \frac{\cc^2}{\sqrt{2}(n+1/2)}(1-|x|)
\;\;\hbox{and}\;\;
|A-1|\le \frac{\cc^2}{\sqrt{3}(n+1/2)}.
\end{equation}
Using \eqref{LP-est-2} and \eqref{est-PSWF} we obtain
\begin{align*}
|\psi_n(x)-\bar{P}_n(x)|\le |\psi_n(x)-A\bar{P}_n(x)| + |A-1||\bar{P}_n(x)|
\le \frac{c'}{n} + \frac{c''}{(n+1)(1-x^2)^{1/4}},
\end{align*}
which implies \eqref{eigen-f-2}.
Inequalities \eqref{norm-psi}, \eqref{est-psi} follow readily from
\eqref{LP-est-1} -- \eqref{eigen-f-2}.
\end{proof}

\subsection{Gaussian bounds for the heat kernel associated to the PSWFs}\label{subsec:PSWF-G}

We now come to one of our principle results on PSWFs of order zero.

\begin{thm}\label{thm:main}
The PSWF heat kernel $p_t(x, y)$ from $\eqref{def-heat-k}$ has the following
Gaussian space localization:
There exist constants $c_1, c_2, c_3, c_4>0$ such that
for any $x, y\in [-1, 1]$ and $t>0$
\begin{equation}\label{Gauss}
\frac{c_1e^{-t\cc^2}\exp\big(- \frac{\rho(x,y)^2}{c_2t}\big)}{\big[V(x, \sqrt t)V(y, \sqrt t)\big]^{1/2}}
\le  p_t(x,y)
\le \frac{c_3\exp\big(- \frac{\rho(x,y)^2}{c_4t}\big)}{\big[V(x, \sqrt t)V(y, \sqrt t)\big]^{1/2}}.
\end{equation}
\end{thm}

\begin{rem}\label{rem:main}
$(a)$
Note that estimates \eqref{Gauss} are valid with
the same constants $c_1$, $c_2$, $c_3$, $c_4$ as in \eqref{Gauss-L0}.
In particular, the upper bound in \eqref{Gauss} is independent of $\cc$.

Because $\rho(x, y):=|\arccos x-\arccos y|$ we have $\rho(x, y) \le \pi$
and hence estimates \eqref{Gauss}
are interesting only in the range $0<t \le\pi$.
Therefore, the factor $e^{-t\cc^2}$ in \eqref{Gauss} is not important;
it can be replaced by $e^{-\pi\cc^2}$ when $0<t \le\pi$.

$(b)$
The {\bf main difference} between the kernels $\KK_{\exp (-tL_\cc)}(x,y)$ and $\KK_{\exp (-tL_0)}(x,y)$
is that $\KK_{\exp (-tL_0)}(x,y)$ obeys the Markov property, while $\KK_{\exp(-tL_\cc)}(x,y)$ does not satisfy it,
namely,
\begin{equation}\label{markov}
\int_{-1}^1 \KK_{\exp (-tL_0)}(x,y)dy \equiv 1,
\quad\hbox{while}\quad
\int_{-1}^1 \KK_{\exp (-tL_\cc)}(x,y)dy \not\equiv 1.
\end{equation}
However, as will be seen below this deficiency of $\KK_{\exp (-tL_\cc)}(x,y)$
is not an obstacle in developing Besov and Triebel-Lizorkin spaces and for all related issues.

\end{rem}

\begin{proof}[Proof of Theorem~\ref{thm:main}]
For this proof we will use our general result from Theorem~\ref{thm:abstract}.
We consider the operator $L_\cc$ as a perturbation of $L_0$, that is,
\begin{equation}\label{def-L-1}
L_\cc f(x):= L_0f(x)+ V(x)f(x),
\quad \quad D(L_\cc) := C^2[-1, 1],
\end{equation}
where $V(x):= \cc^2 x^2\ge 0$.
In the current setting we let
$X:=[-1, 1]$ and $d\mu(x)=dx$
and consider the operators $Z:=L_0$ and $Y:= L_\cc = L_0+V$.

The associated to $L_0$ quadratic form $\cE$ is defined by
\begin{equation}\label{Efg-int}
\cE(f, g):= \int_{\BB^d} L_0 f(x) g(x) dx
= \int_{-1}^1 f'(x)g'(x)dx,
\quad f, g\in D(L_0).
\end{equation}

Our next step is to show that $\overline{\cE}$ is a Dirichlet form.
To this end it suffices to show that (see \cite{Fukushima})
for every $\eps>0$ there exists a function
$\Phi_\eps: \RR\mapsto [-\eps, 1+\eps]$ such that
$\Phi_\eps$ is non-decreasing,
$\Phi_\eps(t)=t$ for $t\in [0, 1]$,
$0\le \Phi_\eps(t') - \Phi_\eps(t) \le t'-t$ if $t<t'$,
and
\begin{equation}\label{BD-cond-1}
f\in D(L_0)
\quad \Longrightarrow\quad
\Phi_\eps(f) \in D(L_0),
\quad
\cE\big(\Phi_\eps(f), \Phi_\eps(f)\big) \le \cE(f, f).
\end{equation}
Let $\Phi_\eps\in C^\infty(\RR)$ have the properties:
$-\eps\le \Phi_\eps\le 1+\eps$, $0\le \Phi_\eps'\le 1$ and $\Phi_\eps(t)=t$, $t\in [0, 1]$,
for some $\eps>0$.
Assume $f\in D(L_0)= C^2[-1, 1]$.
Clearly, $\Phi_\eps(f) \in D(L_0)$ and using that $0\le \Phi_\eps'\le 1$ we obtain
\begin{align*}
\cE(\Phi_\eps(f), \Phi_\eps(f))
&= \int_{-1}^1 \big|[\Phi_\eps(f(x))]'\big|^2 dx
\\
&\le \int_{-1}^1 \big|f'(x)|^2 dx
=\cE(f, f),
\end{align*}
implying \eqref{BD-cond-1}.
Hence, the closure $\overline{\cE}$ of the quadratic form $\cE$ is a Dirichlet form.

Therefore, the theory from Section~\ref{sec:abstract} can be applied.
Invoking Theorem~\ref{thm:abstract} we conclude that
\begin{equation*}
e^{-t\|V\|_\infty}\KK_{\exp(-tL_0)}(x,y) \le p_t(x, y) \le \KK_{\exp(-tL_0)}(x,y).
\end{equation*}
These inequalities and \eqref{Gauss-L0} yield \eqref{Gauss}.
\end{proof}

\subsection{H\"{o}lder continuity of the PSWF heat kernel}\label{subsec:Holder-PSWF}

As is shown in \cite[Section~7]{CKP} $(\overline{\cE}, D(\overline{\cE}))$
with $\overline{\cE}$ and $D(\overline{\cE})$ from above
is a regular strictly local Dirichlet space (see \S \ref{subsec:Holder})
and the respective local scale-invariant Poincar\'{e} inequality is satisfied.
This along with the doubling property \eqref{doubling} yield the Gaussian bounds \eqref{Gauss-L0}
of the kernel $\KK_{\exp(-tL_0)}(x, y)$.
As was indicated in \S \ref{subsec:Holder} inequalities \eqref{Gauss-L0} imply
the H\"{o}lder continuity of $\KK_{\exp(-tL_0)}(x, y)$:
There exists a constant $0<\alpha\le 1$ such that
\begin{equation}\label{Le-Holder}
|\KK_{\exp(-tL_0)}(x, y)-\KK_{\exp(-tL_0)}(x', y)|
\le c_5\Big(\frac{\rho(x, x')}{\sqrt{t}}\Big)^\alpha
\frac{\exp\big(- \frac{\rho(x,y)^2}{c_6t}\big)}{\big[V(x, \sqrt t)V(y, \sqrt t)\big]^{1/2}},
\end{equation}
for all $x, x', y \in [-1, 1]$ and $t>0$, whenever $\rho(x, x')\le \sqrt{t}$.

In turn, applying Theorem~\ref{thm:Holder}, estimates \eqref{Gauss} and \eqref{Le-Holder} imply
that the heat kernel $p_t(x, y)=\KK_{\exp(-tL_\cc)}(x, y)$ is H\"{o}lder continuous as well:
\begin{equation}\label{PSWF-Holder}
|p_t(x, y)-p_t(x', y)|
\le c_7\Big(\frac{\rho(x, x')}{\sqrt{t}}\Big)^\alpha
\frac{\exp\big(- \frac{\rho(x,y)^2}{c_8t}\big)}{\big[V(x, \sqrt t)V(y, \sqrt t)\big]^{1/2}},
\end{equation}
if $\rho(x, x')\le \sqrt{t}$, $0<t\le \pi$.

\begin{rem}\label{rem:Holder-PSWF}
It should be pointed out that the Dirichlet space generated by the operator $L_\cc$ is not strictly local
and, therefore, the H\"{o}lder continuity of $p_t(x, y)$ in $\eqref{PSWF-Holder}$ does not follow by applying
the argument that we used in establishing $\eqref{Le-Holder}$.
\end{rem}

We next establish the following more precise

\begin{thm}\label{thm:Holder-PSWF}
Estimates \eqref{Le-Holder} and \eqref{PSWF-Holder} are valid with $\alpha =1$.
\end{thm}

\begin{proof}
In light of Theorem~\ref{thm:Holder}, we only have to show that
estimate \eqref{Le-Holder} holds with $\gamma =1$.
%To prove this we will proceed similarly as in the proofs of Theorems~\ref{them:ker-D-L} and \ref{thm:ker-lip} below.
We will use the abbreviated notation
$\QQ_t(x, y):=\KK_{\exp(-tL_0)}(x, y)$.

From \eqref{P-eigen-f} and \eqref{heat-ker-L0} it readily follows that
\begin{equation}\label{L-par-p}
L_{0,x} \QQ_t(x,y) = \sum_{n\ge 0} n(n+1)e^{-tn(n+1)} \bar{P}_n(x)\bar{P}_n(y)
= -\partial_t \QQ_t(x, y),
\end{equation}
with $L_{0,x} \QQ_t(x,y)$ denoting the action of $L_0$ on $\QQ_t(x, y)$ as a function of $x$.
It is well known that \eqref{doubling} and \eqref{Gauss-L0} imply
\begin{equation}\label{est-par-p}
|\partial_t \QQ_t(x, y)| \le \frac{c''t^{-1}}{V(y,\sqrt{t})} \exp\Big(- \frac{\rho(x,y)^2}{c't}\Big),
\end{equation}
see e.g. \cite[Theorem 2.32]{GS}, \cite{HS, SC}.
Now, \eqref{def:L0}, \eqref{L-par-p}, and \eqref{est-par-p} yield
\begin{equation}\label{par-p1}
(1-x^2)|\partial_x \QQ_t(x,y)|
\le  \frac{ct^{-1}}{V(y,\sqrt{t})} \int_x^1\exp\Big(- \frac{\rho(u,y)^2}{c't}\Big)du.
\end{equation}
Replacing the integral $\int_x^1$ by $\int_{-1}^x$ we obtain
\begin{equation}\label{par-p2}
(1-x^2)|\partial_x \QQ_t(x,y)|
\le  \frac{ct^{-1}}{V(y,\sqrt{t})} \int_{-1}^x\exp\Big(- \frac{\rho(u,y)^2}{c't}\Big)du.
\end{equation}

We claim that
\begin{equation}\label{dir-Qt}
|\partial_x \QQ_t(x,y)|
\le \frac{ct^{-1/2}(1-x^2)^{-1/2}}{V(y,\sqrt{t})} \exp\Big(- \frac{\rho(x,y)^2}{ct}\Big).
\end{equation}
We will use the following change of variables:
$u=\cos s$,
$y=\cos\phi$, and
$x=\cos \theta$.
Three cases present themselves here depending on the location of $x,y$.

{\em Case 1:} $0\le x\le 1$ and $-1\le y\le x$.
We have
\begin{align*}
&\int_x^1\exp\Big(- \frac{\rho(u,y)^2}{c't}\Big)du
= \int_0^\theta \exp\Big(- \frac{(\phi-s)^2}{c't}\Big)\sin s\, ds
\\
&\le \sin \theta\int_0^\theta \exp\Big(- \frac{1}{c'}\Big(\frac{\phi}{\sqrt{t}}-\frac{s}{\sqrt{t}}\Big)^2\Big) ds
= c\sqrt{t}\sin \theta \int_0^{\theta/\sqrt{t}} \exp\Big(- \frac{1}{c'}\Big(\frac{\phi}{\sqrt{t}}-z\Big)^2\Big) dz
\\
& \le c\sqrt{t}\sin \theta \exp\Big(- \frac{(\phi-\theta)^2}{ct}\Big)
= ct^{1/2}(1-x^2)^{1/2} \exp\Big(- \frac{\rho(x,y)^2}{ct}\Big),
\end{align*}
where we used the well known inequality
$
\int_a^\infty e^{-x^2} dx \le e^{-a^2}, \; a>0.
$
%$$
%\int_a^\infty \exp \Big(-\frac{K x^2}{2}\Big)dx
%\le K^{-1/2}\exp \Big(-\frac{Ka^2}{2}\Big)\Big((\pi/2)^{1/2} \wedge a^{-1}K^{-1/2}\Big), \quad K, a>0.
%$$
The above estimates and \eqref{par-p1} yield \eqref{dir-Qt}.

{\em Case 2:} $0\le x\le 1$ and $x\le y\le 1$ and $t > 1-x^2$.
In this case, we have $\rho(x,y)\le \theta \le c\sqrt{1-x^2} \le c\sqrt{t}$
and hence
$$
\int_x^1\exp\Big(- \frac{\rho(u,y)^2}{c't}\Big)du
\le c(1-x) \le ct^{1/2}(1-x^2)^{1/2}\exp\Big(- \frac{\rho(x,y)^2}{ct}\Big).
$$
This coupled with \eqref{par-p1} implies \eqref{dir-Qt}.

{\em Case 3:} $0\le x\le 1$ and $x\le y\le 1$ and $t\le 1-x^2$.
We have
\begin{align*}
\int_{-1}^x\exp\Big(- \frac{\rho(u,y)^2}{c't}\Big)du
&= \int_\theta^\pi \exp\Big(- \frac{(s-\phi)^2}{c't}\Big)\sin s\, ds
\\
& \le \int_\theta^\pi \exp\Big(- \frac{(s-\phi)^2}{c't}\Big)(\sin \theta +s-\phi) ds.
%\le \int_\theta^\pi \exp\Big(- \frac{(s-\phi)^2}{c't}\Big)(\sin \theta +s-\phi) ds
\end{align*}
Here we used that
$\sin s \le \sin\theta + s-\theta \le \sin\theta + s-\phi$
due to $\phi \le \theta\le s\le \pi$.
Now, using that $\sin\theta = (1-x^2)^{1/2}$ and $\sqrt{t} \le (1-x^2)^{1/2}$ we obtain
\begin{align*}
&\int_{-1}^x\exp\Big(- \frac{\rho(u,y)^2}{c't}\Big)du
\le \int_\theta^\pi \exp\Big(- \frac{(s-\phi)^2}{c't}\Big)\Big(\sin \theta + \sqrt{t}\frac{s-\phi}{\sqrt{t}}\Big) ds
\\
&\le c(1-x^2)^{1/2}\int_\theta^\pi \exp\Big(- \frac{(s-\phi)^2}{2c't}\Big) ds
\le ct^{1/2}(1-x^2)^{1/2}\exp\Big(- \frac{\rho(x,y)^2}{ct}\Big).
\end{align*}
We estimated the last integral above just as the integral in Case 1.
The above estimates and \eqref{par-p2} yield \eqref{dir-Qt}.
The proof of \eqref{dir-Qt} is complete.

%%%%%%%%%%%%%%

\smallskip

We next use \eqref{dir-Qt} to prove \eqref{Le-Holder} with $\alpha=1$.
Two cases are to be considered here.

{\em Case A:} $x, x', y\in[-1,1]$, $x\le x'$, $\rho(x, x')\le \sqrt{t}$, and $\sqrt{t}\le \frac{1}{2}\rho(x, y)$.
Using \eqref{dir-Qt} we obtain
\begin{align*}
|\QQ_t(x, y)-\QQ_t(x', y)|
&\le \int_x^{x'}|\partial_x \QQ_t(u, y)|du
\\
&\le \frac{ct^{-1/2}}{V(y,\sqrt{t})} \int_x^{x'}\frac{1}{\sqrt{1-u^2}} \exp\Big(- \frac{\rho(u,y)^2}{ct}\Big) du.
\end{align*}
Using the assumptions
$
\rho(x, y) \le \rho(x, u)+\rho(u, y) \le \sqrt{t} +\rho(u, y)
\le \frac{1}{2}\rho(x, y)+ \rho(u, y)
$
and hence
$\rho(u, y) \ge \frac{1}{2}\rho(x, y)$.
Therefore,
\begin{align*}
|\QQ_t(x, y)- & \QQ_t(x', y)|
\\
&\le  \frac{ct^{-1/2}}{V(y, \sqrt{t})}|\arccos x-\arccos x'|\exp\Big(- \frac{\rho(x,y)^2}{ct}\Big)
\\
&=\frac{ct^{-1/2}\rho(x, x')}{V(y, \sqrt{t})}\exp\Big(- \frac{\rho(x,y)^2}{ct}\Big).
%\\
%&\le \frac{ct^{-1/2}\rho(x, x')}{\big[V(y, \sqrt{t})V(y, \delta)\big]^{1/2}}
%\exp\Big(- \frac{\rho(x,y)^2}{ct}\Big),
\end{align*}
Using \eqref{double-3} this implies \eqref{Le-Holder} with $\alpha =1$.

\smallskip

{\em Case B:} $x, x', y\in[-1,1]$, $\rho(x, x')\le \sqrt{t}$, and $\frac{1}{2}\rho(x, y)<\sqrt{t}$.
We again use \eqref{dir-Qt} to obtain
\begin{align*}
&|\QQ_t(x, y)- \QQ_t(x', y)|
\le \int_x^{x'}|\partial_x \QQ_t(u, y)|du
\le  \frac{ct^{-1/2}}{V(y, \sqrt{t})}\int_x^{x'}\frac{du}{\sqrt{1-u^2}}
\\
&= \frac{ct^{-1/2}\rho(x, x')}{V(y, \sqrt{t})}
\le \frac{c't^{-1/2}\rho(x, x')}{V(y, \sqrt{t})}\exp\Big(- \frac{\rho(x,y)^2}{4t}\Big),
\end{align*}
which implies \eqref{Le-Holder} with $\alpha =1$.

The proof of Theorem~\ref{thm:Holder-PSWF} is complete.
\end{proof}

\section{Smooth functional calculus induced by the PSWFs of order zero}\label{sec:func-calc}

In this section we develop some basic elements of a smooth functional calculus
associated to the PSWFs of order zero.
Given a bounded function $F: [0, \infty)\to \CC$ the operator $F(\sqrt{L_\cc})$ is defined by
\begin{equation}\label{def-op-fL}
F(\sqrt{L_\cc})g := \sum_{n=0}^\infty F(\sqrt{\chi_n})\langle g, \psi_n\rangle \psi_n,
\quad g\in L^2[-1, 1].
\end{equation}
Then the kernel $\KK_{F(\sqrt{L_\cc})}(x, y)$ of $F(\sqrt{L_\cc})$ (if it exists)  takes the form
\begin{equation}\label{kernel-fL}
\KK_{F(\sqrt{L_\cc})}(x, y)= \sum_{n=0}^\infty F(\sqrt{\chi_n})\psi_n(x)\psi_n(y),
\;\; x, y\in [-1, 1]
\;\;\; (\overline{\psi_n(x)}=\psi_n(x)).
\end{equation}

Our goal is to establish localization results for the kernel $\KK_{F(\sqrt{L_\cc})}(x, y)$
in the case when the function $F$ is smooth and with compact support or rapid decay.

\subsection{Main localization results}\label{subsec:smooth-f-c}

\begin{thm}\label{thm:ker-local}
Let $F\in C^k(\RR)$, $k\ge 3$,  $F$ even, real-valued, and $\supp F\subset [-R, R]$ for some $R\ge 1$.
Then $F(\delta\sqrt{L_\cc})$, $\delta>0$, is an integral operator with kernel $\KK_{F(\delta\sqrt{L_\cc})}(x,y)$
satisfying
\begin{equation}\label{ker-local}
|\KK_{F(\delta\sqrt{L_\cc})}(x,y)| \le \frac{c_1\Big(1+\frac{\rho(x, y)}{\delta}\Big)^{-k}}{[V(x, \delta)V(y, \delta)]^{1/2}},
\quad x, y\in [-1, 1],\; \delta >0,
\end{equation}
where
$c_1=c_2R^2(\|F\|_\infty+R^k\|F^{(k)}\|_\infty)$
with the constant $c_2$ depending only on $k$.
\end{thm}

As a consequence of Theorem~\ref{thm:ker-local} we obtain the following localization result (see the proof of Theorem~3.4 in \cite{KP}),
where we use the Schwartz class $\cS(\RR)$, defined by
\begin{equation}\label{def-S}
\cS(\RR):= \big\{\phi\in C^\infty(\RR): \sup_{x\in\RR}|x^k\phi^{(n)}(x)|<\infty, \;\;\forall k, n\in \NN_0\big\}.
\end{equation}

\begin{thm}\label{thm:ker-local-S}
If the function $F$ is even, real-valued, and $F\in \cS(\RR)$, then for any $\sigma>0$
\begin{equation}\label{ker-local-S}
|\KK_{F(\delta\sqrt{L_\cc})}(x,y)|
\le \frac{c_\sigma\Big(1+\frac{\rho(x, y)}{\delta}\Big)^{-\sigma}}{\big[V(x, \delta)V(y, \delta)\big]^{1/2}},
\quad x, y\in [-1, 1],\; \delta >0.
\end{equation}
\end{thm}

\begin{rem}\label{rem:denominator}
Because of inequality \eqref{double-3}, the denominator $\big[V(x, \delta)V(y, \delta)\big]^{1/2}$
in \eqref{ker-local-S} can be replaced by $V(x, \delta)$ or $V(y, \delta)$
$($with a different value of $\sigma$, which does not change the statement$)$.
The same is valid for the next result.
\end{rem}

\begin{rem}\label{rem:sub-exp}
The space localization of the kernels $\KK_{F(\delta\sqrt{L_\cc})}(x,y)$ can be improved to sub-exponential.
Just as in \cite[Theorem~3.6]{KP} it can be shown that
for any $0<\varepsilon \le 1$ there exists an even ``cut-off" function $F$ with ``small derivatives"
such that
\begin{equation}\label{exp-loc}
|\KK_{F(\delta\sqrt{L_\cc})}(x,y)|
\le \frac{c_1\exp\Big(-c_2\big(\frac{\rho(x, y)}{\delta}\big)^{1-\varepsilon}\Big)}{[V(x, \delta)V(x, \delta)]^{1/2}},
\quad \delta>0.
\end{equation}
\end{rem}

\subsection{Proof of the main localization estimate (Theorem~\ref{thm:ker-local})}

The proof of Theorem~\ref{thm:ker-local} follows in the foot steps of the proof of Theorem~3.1 in \cite{KP}.
For this proof we need several ingredients.
First, we need an estimate for $|\KK_{F(\sqrt{L_\cc})}(x, y)|$
when $F$ is ``rough" but compactly supported.

\begin{prop}\label{prop:rough-oper}
$(a)$
Let $F:[0, \infty)\to \CC$ be bounded and $\supp F\subset [0, \tau]$ for some $\tau\ge 1$.
Then $F(\sqrt{L_\cc})$ is an integral operator whose kernel $\KK_{F(\sqrt{L_\cc})}(x, y)$ is continuous
and
\begin{equation}\label{rough-ker-1}
|\KK_{F(\sqrt{L_\cc})}(x, y)| \le \frac{c\|F\|_\infty}{V(x, \tau^{-1})^{1/2}V(y, \tau^{-1})^{1/2}},
\quad x, y\in [-1, 1].
\end{equation}

$(b)$
Let $F:[0, \infty)\to \CC$ be bounded and $|F(\lambda)| \le c\lambda^{-\sigma}$, $\lambda\ge 1$, for some $\sigma>2$.
Then $F(\sqrt{L_\cc})$ is an integral operator whose kernel $\KK_{F(\sqrt{L_\cc})}(x, y)$ is continuous on $[-1, 1]^2$.
\end{prop}

\begin{proof}
(a) From \eqref{norm-psi} and \eqref{est-psi} along with \eqref{mes-B-inv}
it readily follows that for $\tau \ge 1$
\begin{equation*} %\label{rough-ker-22}
\sum_{0\le n\le \tau} \psi_n^2(x)
\le c\tau\min\Big\{\frac{1}{(1-x^2)^{1/2}}, \frac{1}{\tau^{-1}}\Big\}
\le \frac{c}{V(x, \tau^{-1})}
\end{equation*}
and using also \eqref{eigen-v} we obtain
\begin{equation}\label{rough-ker-23}
\sum_{n\ge 0}|F(\sqrt{\chi_n})|\psi_n^2(x)
\le \|F\|_\infty\sum_{0\le n\le \tau}\psi_n^2(x)
\le \frac{c\|F\|_\infty}{V(x, \tau^{-1})}.
\end{equation}
From this applying the Cauchy-Schwarz inequality we obtain
\begin{align*}
|\KK_{F(\sqrt{L_\cc})}(x, y)| &\le \Big(\sum_{n\ge 0}|F(\sqrt{\chi_n})|\psi_n^2(x)\Big)^{1/2}
\Big(\sum_{n\ge 0}|F(\sqrt{\chi_n})|\psi_n^2(y)\Big)^{1/2}
\\
&\le \frac{c\|F\|_\infty}{V(x, \tau^{-1})^{1/2}V(y, \tau^{-1})^{1/2}},
\end{align*}
which confirms \eqref{rough-ker-1}.

The kernel $\KK_{F(\sqrt{L_\cc})}(x, y)$ is continuous because the PSWFs $\{\psi_n\}$ are continuous.

\smallskip

(b) Under the assumptions of part (b) we have
\begin{align*}
|\KK_{F(\sqrt{L_\cc})}(x, y)|
&\le  \sum_{n=0}^\infty |F(\sqrt{\chi_n})||\psi_n(x)||\psi_n(y)|
\\
&\le c \sum_{n=0}^\infty \chi_n^{-\sigma/2}(n+1)
\le c \sum_{n=1}^\infty n^{-\sigma+1} < \infty,
\end{align*}
were we used \eqref{eigen-v} and \eqref{norm-psi}.
Consequently, the above series converges uniformly and since the functions $\{\psi_n\}$ are continuous
the kernel $\KK_{F(\sqrt{L_\cc})}(x, y)$ is also continuous.
\end{proof}

An important component in the proof of Theorem~\ref{thm:ker-local} is the following

\subsection*{Davies-Gaffney estimate}

We say that
the PSWF semigroup $\Psi(z):=\exp(-zL_\cc)$, $z\in \CC_+$, %in the setting described in the introduction
satisfies the Davies-Gaffney estimate if
\begin{equation}\label{D-G-Lc}
|\langle \Psi(t)f_1, f_2\rangle| \le \exp\Big(-\frac{\hat{c}r^2}{t}\Big)\|f_1\|_2 \|f_2\|_2
\end{equation}
for all $t>0$, $U_j\subset [-1, 1]$, $U_j$ open, $f_j\in L^2[-1, 1]$,
$\supp f_j\subset U_j$, $j=1, 2$, and $r=\rho(U_1, U_2)$,
where $\hat{c}>0$ is a constant.

As a consequence of the Gaussian upper bound for the heat kernel from \eqref{gauss-ball-Lc}
we get the following

\begin{prop}\label{prop:D-G}
The Davies-Gaffney estimate holds for the PSWF semigroup $\Psi(z):=\exp(-zL_\cc)$
with constant $\hat{c}=c_4^{-1}$, where $c_4$ is the constant from \eqref{Gauss}.
\end{prop}

The proof of this proposition is a repetition of the proof of Proposition~2.7 in \cite{KP};
we omit it.

\subsection*{Finite speed propagation property of the wave equation}

\begin{prop}\label{prop:f-speed}
In our setting the Davies-Gaffney estimate \eqref{D-G-Lc} is equivalent to the so called
finite speed propagation property:
\begin{equation}\label{def-f-speed-2}
\langle \cos(t\sqrt{L_\cc})f_1, f_2 \rangle =0,
\quad 0<\tilde{c}t<r, \quad \tilde{c} := \frac{1}{2\sqrt{\hat{c}}},
\end{equation}
for all open sets $U_j\subset [-1, 1]$, $f_j\in L^2[-1, 1]$, $\supp f_j\subset U_j$, $j=1, 2$, $r:=\rho(U_1, U_2)$.
\end{prop}

For the proof of this proposition, see the proof of Theorem~3.4 in \cite{CS}.

In turn Proposition~\ref{prop:f-speed} yields the following

\begin{prop}\label{prop:finite-speed}
Let $F$, defined on $\RR$, be even, $\supp \widehat{F}\subset [-A, A]$ for some $A>0$, and $\widehat{F}\in W_1^m$ for some $m>2$,
i.e. $\|\widehat{F}^{(m)}\|_1<\infty$.
Here $\widehat{F}(\xi):=\int_\RR F(u)e^{-iu\xi}du$.
Then for $\delta>0$ and $x, y\in [-1, 1]$ we have
\begin{equation}\label{f-s-ker}
\KK_{F(\delta\sqrt{L_\cc})}(x, y) = 0
\quad\hbox{if}\quad
\rho(x, y)>\delta A.
\end{equation}
\end{prop}

The proof of this proposition is carried out just like the proof of \cite[Proposition~2.8]{KP}
with the only difference that we use here Proposition~\ref{prop:rough-oper} (b) to conclude that
the kernel $\KK_{F(\delta\sqrt{L_\cc})}(x, y)$ is continuous.
We omit the further details.

\subsection*{Completion of the proof of Theorem~\ref{thm:ker-local}}

The proof of Theorem~\ref{thm:ker-local}
follows by Proposition~\ref{prop:rough-oper} and Proposition~\ref{prop:finite-speed}.
The proof is identical to the proof of the respective part of \cite[Theorem 3.1]{KP}
and will be omitted.
\qed

\subsection{Lipschitz continuity and smoothness of kernels}\label{subsec:Lip-kernel}

We now study how the differential operators
\begin{equation}\label{operators}
\sD_x :=(1-x^2)^{1/2}\partial_x
\quad\hbox{and}\quad
L_{0,x}:= -\partial_x[(1-x^2)\partial_x]
\end{equation}
affect the space localization of kernels $\KK_{F(\delta\sqrt{L_\cc})}(x, y)$.
The operators $\sD_y$ and $L_{0, y}$ are defined similarly.

\begin{thm}\label{them:ker-D-L}
Let $F\in \cS(\RR)$ be even and real-valued.
Then for any $\sigma>0$ there exists a constant $c_\sigma>0$ such that
\begin{equation}\label{ker-D}
|\sD_x \KK_{F(\delta\sqrt{L_\cc})}(x,y)|
\le \frac{c_\sigma \delta^{-1}\Big(1+\frac{\rho(x, y)}{\delta}\Big)^{-\sigma}}{[V(x, \delta)V(y, \delta)]^{1/2}},
\quad x, y\in [-1, 1],\; 0<\delta\le 1,
\end{equation}
and
\begin{equation}\label{ker-L}
|L_{0,x} \KK_{F(\delta\sqrt{L_\cc})}(x,y)|
\le \frac{c_\sigma \delta^{-2}\Big(1+\frac{\rho(x, y)}{\delta}\Big)^{-\sigma}}{[V(x, \delta)V(y, \delta)]^{1/2}},
\quad x, y\in [-1, 1],\; 0<\delta\le 1.
\end{equation}
Because of the symmetry of $\KK_{F(\delta\sqrt{L_\cc})}(x,y)$
the above inequalities are valid with the operators
$\sD_x$, $L_{0, x}$ replaced by $\sD_y$, $L_{0, y}$.
\end{thm}

\begin{proof}
To prove \eqref{ker-D} we observe that
from \eqref{def:Lc}, \eqref{eigenf-L}, and \eqref{def:L0} it follows that
\begin{equation}\label{lip-0}
L_0\psi_n(x) = L_\cc\psi_n(x) - \cc^2x^2\psi_n(x) = \chi_n\psi_n(x) - \cc^2x^2\psi_n(x).
\end{equation}
Using this, \eqref{ker-fLc}, and the rapid decay of $\{F(\delta\sqrt{\chi_n})\}$ we get
\begin{align*}
L_{0, x}&\KK_{F(\delta\sqrt{L_\cc})}(x, y)= \sum_{n=0}^\infty F(\delta\sqrt{\chi_n})L_{0, x}\psi_n(x)\psi_n(y)
\\
&= \sum_{n=0}^\infty \chi_n F(\delta\sqrt{\chi_n})\psi_n(x)\psi_n(y)
- \cc^2 x^2\sum_{n=0}^\infty F(\delta\sqrt{\chi_n})\psi_n(x)\psi_n(y).
\end{align*}
Hence
\begin{align*}
\partial_x\big[(1-x^2)\partial_x\KK_{F(\delta\sqrt{L_\cc})}(x, y)\big]
= - \KK_{L_\cc F(\delta\sqrt{L_\cc})}(x, y) + \cc^2x^2\KK_{F(\delta\sqrt{L_\cc})}(x, y)
\end{align*}
and integrating with respect to $x$ we arrive at
\begin{align}\label{dx-int}
(1-x^2)\partial_x \KK_{F(\delta\sqrt{L_\cc})}(x,y)
= \int_x^1 \big[\KK_{L_\cc F(\delta\sqrt{L_\cc})}(u, y) - \cc^2u^2 \KK_{F(\delta\sqrt{L_\cc})}(u,y)\big] du.
\end{align}
Consider the function $H(\lambda):=\lambda^2F(\lambda)$.
Clearly, $H(\delta \sqrt{L_\cc})=\delta^2 L_\cc F(\delta \sqrt{L_\cc})$.
On the other hand, $H\in\cS(\RR)$ and $H$ is even and real-valued.
Then by Theorem~\ref{thm:ker-local-S}, applied to $H$, we conclude that for any $\sigma>0$
\begin{equation*}
\big|\KK_{L_c f(\delta\sqrt{L_\cc})}(x,y)\big|
\le \frac{c_\sigma\delta^{-2}\Big(1+\frac{\rho(x, y)}{\delta}\Big)^{-\sigma}}{\big[V(x, \delta)V(y, \delta)\big]^{1/2}}
\le \frac{c\delta^{-2}}{V(y, \delta)}
\Big(1+\frac{\rho(x, y)}{\delta}\Big)^{-\sigma+1},
%\quad x, y\in [-1, 1],\; \delta >0.
\end{equation*}
where for the last inequality we used \eqref{double-3}.
Again by Theorem~\ref{thm:ker-local-S}, it follows that
\begin{equation*}
|\KK_{F(\delta\sqrt{L_\cc})}(x,y)|
\le \frac{c}{V(y, \delta)}
\Big(1+\frac{\rho(x, y)}{\delta}\Big)^{-\sigma}.
\end{equation*}
Putting the above together we obtain
\begin{equation}\label{first-est}
(1-x^2)^{1/2}|\partial_x \KK_{F(\delta\sqrt{L_\cc})}(x,y)|
\le \frac{c\delta^{-2}}{V(y, \delta)}\frac{1}{(1-x^2)^{1/2}}
\int_x^1 \Big(1+\frac{\rho(u, y)}{\delta}\Big)^{-\sigma+1}du.
\end{equation}
Replacing the integral $\int_x^1$ above by $\int_{-1}^x$ we get
\begin{equation}\label{second-est}
(1-x^2)^{1/2}|\partial_x \KK_{F(\delta\sqrt{L_\cc})}(x,y)|
\le \frac{c\delta^{-2}}{V(y, \delta)}\frac{1}{(1-x^2)^{1/2}}
\int_{-1}^x \Big(1+\frac{\rho(u, y)}{\delta}\Big)^{-\sigma+1}du.
\end{equation}

In light of \eqref{double-3} to prove \eqref{ker-D} it suffices to show that for any $\sigma>0$
\begin{equation}\label{ker-D-2}
(1-x^2)^{1/2}|\partial_x \KK_{F(\delta\sqrt{L_\cc})}(x,y)|
\le \frac{c\delta^{-1}}{V(y, \delta)}\Big(1+\frac{\rho(x, y)}{\delta}\Big)^{-\sigma+3}
\end{equation}
with the constant $c>0$ depending on $\sigma$.
On the other hand, from \eqref{first-est} and \eqref{second-est} it follows that
to prove \eqref{ker-D-2} it is sufficient to prove this estimate
\begin{align}\label{x-int-1}
\min\Big (\int_x^1 \Big(1+\frac{\rho(u, y)}{\delta}\Big)^{-\sigma+1}du,
&\int_{-1}^x \Big(1+\frac{\rho(u, y)}{\delta}\Big)^{-\sigma+1}du\Big)
\\
&\le c\delta (1-x^2)^{1/2} \Big(1+\frac{\rho(x, y)}{\delta}\Big)^{-\sigma+3}. \nonumber
\end{align}
Moreover, because of the symmetry, we may assume that $x\geq 0$.
In what follows we assume that $\sigma>0$ is sufficiently large.

We will distinguish three cases, where we apply the same change of variables:
$u=\cos \alpha$, $y=\cos \phi$, and $x=\cos \theta$.
\smallskip

{\em Case 1.} Let $0\le x\le 1$ and $-1\le y\le x$.
Then
\begin{align*}
& \int_x^1 \Big(1+\frac{\rho(u, y)}{\delta}\Big)^{-\sigma+1}du
=\int_0^\theta \frac{\sin\alpha}{\Big(1+\frac{\phi-\alpha}{\delta}\Big)^{\sigma-1}}d\alpha\leq \int_0^\theta \frac{\sin\theta}{\Big(1+\frac{\phi-\alpha}{\delta}\Big)^{\sigma-1}}d\alpha
\\
&\quad
\le c\delta \sin \theta\Big(1+\frac{\phi-\theta}{\delta}\Big)^{-\sigma+2}
= c\delta (1-x^2)^{1/2} \Big(1+\frac{\rho(x, y)}{\delta}\Big)^{-\sigma+2},
\end{align*}
which confirms \eqref{x-int-1}.

\smallskip

{\em Case 2.} Let $0\le x\le 1$, $x\le y\le 1$, and $\delta >(1-x^2)^{1/2}$.
The last condition implies that $\theta\leq c\delta$ and
hence $\rho (x,y) \le \theta\le c\delta$ whenever $x\le u \le 1$.
Therefore,
$$
\int_x^1 \Big(1+\frac{\rho(u, y)}{\delta}\Big)^{-\sigma+1}du\leq c(1-x)
\le c (1-x^2)^{1/2}\delta \Big(1+\frac{\rho(x, y)}{\delta}\Big)^{-\sigma+2},
$$
which confirms \eqref{x-int-1} in this case.

\smallskip

{\em Case 3.} Let $0\le x\le 1$, $x\le y\le 1$, and $\delta\leq (1-x^2)^{1/2}$.
We now consider the other integral.
\begin{align*}\label{x-int-2}
 \int_{-1}^x \Big(1+\frac{\rho(u, y)}{\delta}\Big)^{-\sigma+1}du
= \int_\theta^\pi \frac{\sin \alpha}{\Big(1+\frac{\alpha-\phi}{\delta}\Big)^{\sigma-1}} d\alpha
\leq \int_\theta^\pi \frac{\sin \theta +\alpha-\phi}{\Big(1+\frac{\alpha-\phi}{\delta}\Big)^{\sigma-1}} d\alpha,
\end{align*}
where we used the inequality
$\sin \alpha\le \sin \theta + \alpha-\theta\le \sin \theta + \alpha-\phi$,
taking into account that
$\phi\le \theta\le \alpha\le \pi$.
The first term, containing $\sin \theta$, is treated as before.
For the second one, we write
$\alpha-\phi\le \delta(1+\frac{\alpha-\phi}{\delta})$.
Hence,
\begin{align*}
\int_{-1}^x &\Big(1+\frac{\rho(u, y)}{\delta}\Big)^{-\sigma+1}du
\le \sin \theta \int_{\theta}^\pi \Big(1+\frac{\alpha-\phi}{\delta}\Big)^{-\sigma+1} d\alpha
\\
&+\delta\int_{\theta}^\pi  \Big(1+\frac{\alpha-\phi}{\delta}\Big)^{-\sigma+2} d\alpha
\le \frac{c\delta(\sin \theta +\delta) )}{\Big(1+\frac{\theta-\phi}{\delta}\Big)^{\sigma-3}}
=\frac{c\delta((1-x^2)^{1/2} +\delta)}{\Big(1+\frac{\rho(x, y)}{\delta}\Big)^{\sigma-3}},
\end{align*}
which along with the assumption on $\delta$ implies \eqref{x-int-1}.

This completes the proof of estimate \eqref{ker-D-2},
which in turn implies \eqref{ker-D}.

The proof of \eqref{ker-L} is easier. There is no integration.
From \eqref{lip-0} it follows that
\begin{align*}
L_{0,x} \KK_{F(\delta\sqrt{L_\cc})}(x,y)
= -\cc^2x^2 \KK_{F(\delta\sqrt{L_\cc})}(x,y) + \KK_{L_\cc f(\delta\sqrt{L_\cc})}(x, y)
\end{align*}
and estimate \eqref{ker-L} follows immediately applying Theorem~\ref{thm:ker-local-S}.
\end{proof}

We next show the Lipschitz continuity of kernels $\KK_{F(\delta\sqrt{L_\cc})}(x, y)$ for smooth $F$'s.

\begin{thm}\label{thm:ker-lip}
Let $F\in \cS(\RR)$ be even and real-valued.
Then for any $\sigma>0$ there exists a constant $c_\sigma>0$ such that
for all $x,x', y\in [-1, 1]$, $\delta >0$
\begin{equation}\label{ker-lip-11}
|\KK_{F(\delta\sqrt{L_\cc})}(x,y)- \KK_{F(\delta\sqrt{L_\cc})}(x',y)|
\le c_\sigma\frac{\frac{\rho(x, x')}{\delta}\Big(1+\frac{\rho(x, y)}{\delta}\Big)^{-\sigma}}{[V(x, \delta)V(y, \delta)]^{1/2}},
\;\;\hbox{if}\;\; \rho(x, x')\le \delta,
\end{equation}
\end{thm}

\begin{proof}
Two cases present themselves here.

(i) Assume $x, x', y\in[-1,1]$, $x\le x'$, $\rho(x, x')\le \delta$, and $0<\delta\le \frac{1}{2}\rho(x, y)$.
Then using \eqref{ker-D} we get
\begin{align*}
|\KK_{F(\delta\sqrt{L_\cc})}(x, y)-\KK_{F(\delta\sqrt{L_\cc})}(x', y)|
&\le \int_x^{x'}|\partial_x \KK_{F(\delta\sqrt{L_\cc})}(u, y)|du
\\
&\le  \frac{c\delta^{-1}}{V(y, \delta)}\int_x^{x'}\frac{1}{\sqrt{1-u^2}} \Big(1+\frac{\rho(u, y)}{\delta}\Big)^{-\sigma}du.
\end{align*}
However,
$
\rho(x, y) \le \rho(x, u)+\rho(u, y) \le \delta +\rho(u, y)
\le \frac{1}{2}\rho(x, y)+ \rho(u, y)
$
and hence
$\rho(u, y) \ge \frac{1}{2}\rho(x, y)$,
implying
\begin{align*}
|\KK_{F(\delta\sqrt{L_\cc})}(x, y)- & \KK_{F(\delta\sqrt{L_\cc})}(x', y)|
\\
&\le  \frac{c\delta^{-1}}{V(y, \delta)}|\arccos x-\arccos x'|\Big(1+\frac{\rho(x, y)}{\delta}\Big)^{-\sigma}
\\
&=\frac{c\delta^{-1}\rho(x, x')}{V(y, \delta)}\Big(1+\frac{\rho(x, y)}{\delta}\Big)^{-\sigma}
\\
&\le \frac{c\delta^{-1}\rho(x, x')}{\big[V(y, \delta)V(y, \delta)\big]^{1/2}}
\Big(1+\frac{\rho(x, y)}{\delta}\Big)^{-\sigma+1},
\end{align*}
which implies \eqref{ker-lip-11} because the constant $\sigma>0$ can be arbitrarily large.
Above for the last inequality we used \eqref{double-3}.

(ii) Assume $x, x', y\in[-1,1]$, $\rho(x, x')\le \delta$, and $\frac{1}{2}\rho(x, y)<\delta\le 1$.
Then using again \eqref{ker-D} we obtain
\begin{align*}
&|\KK_{F(\delta\sqrt{L_\cc})}(x, y)- \KK_{F(\delta\sqrt{L_\cc})}(x', y)|
\le \int_x^{x'}|\partial_x \KK_{F(\delta\sqrt{L_\cc})}(u, y)|du
\\
&\qquad\le  \frac{c\delta^{-1}}{V(y, \delta)}\int_x^{x'}\frac{du}{\sqrt{1-u^2}}
= \frac{c\delta^{-1}\rho(x, x')}{V(y, \delta)}
\le \frac{c'\delta^{-1}\rho(x, x')}{V(y, \delta)}\Big(1+\frac{\rho(x, y)}{\delta}\Big)^{-\sigma},
\end{align*}
which implies \eqref{ker-lip-11}.
\end{proof}

\section{Besov and Triebel-Lizorkin spaces associated to PSWFs}\label{sec:B-F-spaces}

As elsewhere the operators $L_0$ and $L_\cc$, $\cc>0$, see \eqref{def:L0} and \eqref{def:Lc},
will play a major role.
In this section we introduce Besov and Triebel-Lizorkin spaces associated to the operators $L_0$ and $L_c$.
One of our main goals will be to show that these spaces are the same when $L_0$ is replaced by $L_c$.

\subsection{Distributions}\label{subsec:distributions}

The Besov and Triebel-Lizorkin spaces that we are interested in are in general spaces of distributions.
Here we introduce and discuss distributions that are relevant to our setting.

It is natural to define the set $\cD$ of test functions as the set of all functions $\phi\in C^\infty(-1, 1)$ such that
$\|L_0^m\phi\|_2 < \infty$ for all $m\in \NN_0$.
The topology on $\cD$ is determined by
\begin{equation}\label{def-Pm}
\cP_m(\phi):= \|L_0^m\phi\|_2, \quad m\ge 0
\quad\hbox{or by}\quad
\cP_m^\star(\phi):= \max_{0\le r\le m}\|L_0^r\phi\|_2,
\quad m\ge 0.
\end{equation}

Observe that for $\phi\in\cD$ we have
$
\langle \phi, \bar{P}_n\rangle = [n(n+1)]^{-k}\langle L_0^k\phi, \bar{P}_n\rangle
$
and hence
\begin{equation*}
|\langle \phi, \bar{P}_n\rangle|\le n^{-2k}\|L_0^k\phi\|_2,
\quad \forall k\in\NN_0.
\end{equation*}
On the other hand, by Parseval's identity
\begin{align*}
\|L_0^k\phi\|_2 = \Big(\sum_{n\ge 0} |\langle \phi, \bar{P}_n\rangle|^2[n(n+1)]^{4k}\Big)^{1/2},
\quad k\in\NN.
\end{align*}
Therefore, the topology on $\cD$ is completely determined by the Legendre-Fourier coefficients of $\phi$.
Namely, the topology on $\cD$ can be equivalently defined by
\begin{equation}\label{def-Pm-star}
\cP_m^{\star\star}(\phi):= \sup_{n\ge 0}n^m |\langle \phi, \bar{P}_n\rangle|,
\quad m\ge 0.
\end{equation}

The space $\cD'$ of distributions in our setting is defined as the space of all continuous linear functionals on $\cD$.
The pairing of $f\in\cD'$ and $\phi\in\cD$ will be denoted by
$\langle f, \phi\rangle := f(\overline{\phi})$;
it is consistent with the inner product
$\langle f, g\rangle :=\int_I f\overline{g} dx$ in $L^2(I)$, $I:=[-1, 1]$.
Note that a linear functional $f\in \cD'$ if and only if there exist constants $m\ge 0$ and $c_m>0$ such that
\begin{equation}\label{distr-1}
|\langle f, \phi\rangle| \le c_m\cP_m^\star(\phi), \quad \forall \phi\in \cD.
\end{equation}

We introduce the operator
\begin{equation}\label{def-D}
\fD f(x):=(1-x^2)\frac{df}{dx}(x)
\end{equation}
and will need this lemma:

\begin{lem}\label{lem:L0-D-Lc}
For any $k\in\NN$ and $f\in C^{2k}(-1, 1)$
\begin{align}
L_\cc^k f
&=\sum_{\nu=0}^k \alpha_\nu L_0^\nu f + \sum_{\nu=0}^{k-1} \beta_\nu\fD L_0^\nu f,
\quad\hbox{and} \label{rep-Lc-L0D}
\\
L_0^k f
&=\sum_{\nu=0}^k \gamma_\nu L_\cc^\nu f + \sum_{\nu=0}^{k-1} \eta_\nu \fD L_\cc^\nu f, \label{rep-L0-LcD}
\end{align}
where $\alpha_k=\gamma_k=1$ and
$\alpha_\nu$, $\beta_\nu$, $\gamma_\nu$, and $\eta_\nu$ for $\nu=0, \dots, k-1$
are polynomials of degree $\le 3k$.
\end{lem}

\noindent
{\em Proof.}
Identities \eqref{rep-Lc-L0D}, \eqref{rep-L0-LcD} are easily proved by induction
using the obvious identities:
$L_cf(x) = L_0 f(x)+\cc^2x^2f(x)$, $L_0f(x) = L_\cc f(x)-\cc^2x^2f(x)$,
\begin{align*}
L_0(fg)&= gL_0f+fL_0g - 2g'\fD f,
\quad
L_0\fD f(x)= \fD L_0f(x) +2xL_0f(x),
\\
\fD(fg)&=g\fD f + f\fD g,
\quad
\fD^2f(x)=-(1-x^2)L_0f(x).
\hspace{1.5in} \hbox{$\qed$}
\end{align*}

From Lemma~\ref{lem:L0-D-Lc} it follows that
$\phi\in \cD$ if and only if $\|L_\cc^m \phi\|_2 <\infty$ for all $m\ge 0$.
In turn, this is equivalent to the condition
$\sup_{n\ge 0} n^m|\langle \phi, \psi_n\rangle| <\infty$ for all $m\ge 0$;
the proof is just as in the case of the Legendre polynomials above.
Therefore, in the definition of distributions above the role of $L_0$ can be played by $L_\cc$
and the role of Legendre polynomials $\{\bar{P}_n\}$ can be played by the PSWFs $\{\psi_n\}$.

We will need the following basic convergence result.

\begin{prop}\label{prop:converge}
Assume $\varphi\in C^\infty(\RR)$, $\varphi$ is even, real-valued,
$\supp \varphi \subset [-R, R]$, $R>0$, and $\varphi(0)=1$.
Then for any $f\in \cD'$ $($or $f\in \cD$$)$
\begin{equation}\label{conv-1}
f=\lim_{\delta\to 0} \varphi(\delta\sqrt{L_0})f
\quad\hbox{in $\cD'$ $($or in $\cD$$)$,}
\end{equation}
and also
\begin{equation}\label{conv-2}
f=\lim_{\delta\to 0} \varphi(\delta\sqrt{L_\cc})f
\quad\hbox{in $\cD'$ $($or in $\cD$$)$.}
\end{equation}
\end{prop}

The proof of this proposition is similar to the proof of Proposition~5.5 in \cite{KP}.
We omit it.

\subsection{Maximal operators}\label{subsec:max-oper}

We will utilize the following version of the Hardy-Littlewood maximal operator:
\begin{equation}\label{def-max-oper}
\cM_tf(x):=\sup_{r>0}\Big(\frac{1}{V(x,r)}\int_{B(x,r)}|f(y)|^tdy\Big)^{1/t},
\quad x\in [-1, 1], \; t>0.
\end{equation}
The doubling property \eqref{doubling} implies that
the Fefferman-Stein vector-valued maximal inequality holds:
If~$0<p<\infty$, $0<q\le \infty$, and $0<t<\min\{p, q\}$, then for any sequence of functions $\{f_\nu\}$ on $[-1, 1]$
\begin{equation}\label{vector-max}
\Big\|\sum_\nu \Big(|\cM_tf_\nu(\cdot)|^q\Big)^{1/q}\Big\|_{L^p}
\le c\Big\|\Big(\sum_\nu|f_\nu(\cdot)|^q\Big)^{1/q}\Big\|_{L^p}.
\end{equation}

We will also need estimates on Peetre-type maximal operators.
Denote by $\Pi_N$ the set of all univariate algebraic polynomials of degree $\le N$
and also define
\begin{equation}\label{def-Sigma-N}
\Sigma_N^\cc:=\Big\{g: g=\sum_{n=0}^Na_n\psi_n,\; a_n\in\RR \Big\}.
\end{equation}
Thus $\{\Pi_N\}$ and $\{\Sigma_N^\cc\}$ are the spectral spaces associated
with the operators $L_0$ and~$L_\cc$, respectively.

\begin{lem}\label{lem:peetre}
If $g\in \Pi_N$ or $g\in\Sigma_N^\cc$, $N\ge 0$,
then
\begin{equation}\label{peetre}
\sup_{y\in [-1, 1]}\frac{|g(y)|}{(1+N\rho(x, y))^{2/t}}
\le c \cM_t g(x),
\quad x\in [-1, 1], \; t>0.
\end{equation}
\end{lem}

For the proof of this lemma, see Lemma~6.4 in \cite{KP}.

\subsection{Besov spaces associated to $L_0$ and $L_\cc$}\label{subsec:Besov-L0}

Following ideas from \cite{KP} we introduce two kinds of Besov spaces associated with $L_0$:

(a) Classical Besov spaces $B^s_{pq}=B^s_{pq}(L_0)$, and

(b) Nonclassical Besov spaces $\widetilde{B}^s_{pq}=\widetilde{B}^s_{pq}(L_0)$.

Let $\varphi_0, \varphi \in C^\infty(\RR)$ be even real-valued functions such that
\begin{align}
&\supp \varphi_0 \subset [0, 2],
\quad
|\varphi_0(\lambda)| \ge c>0 \;\; \hbox{for}\;\; \lambda\in [0, 2^{3/4}], \label{phi-0}
\\
&\supp \varphi \subset [1/2, 2],
\quad
|\varphi(\lambda)| \ge c>0 \;\; \hbox{for}\;\; \lambda\in [2^{-3/4}, 2^{3/4}]. \label{phi}
\end{align}
Hence $|\varphi(\lambda)| + \sum_{j\ge 1}|\varphi(2^{-j}\lambda)\ge c>0$ for $\lambda \in \RR_+$.
Set $\varphi_j(\lambda):=\varphi(2^{-j}\lambda)$, $j\ge 1$.

\begin{defn}\label{def:Besov-0}
Let $s\in\RR$ and $0< p, q\le \infty$.

$(a)$ The Besov space $B^s_{pq}=B^s_{pq}(L_0)$ is defined as the set of all $f\in \cD'$ such that
\begin{equation}\label{Besov-1}
\|f\|_{B^s_{pq}} := \Big(\sum_{j\ge 0}\big(2^{sj}\|\varphi_j(\sqrt{L_0})f\|_p\big)^q\Big)^{1/q} < \infty.
\end{equation}

$(b)$ The Besov space $\widetilde{B}^s_{pq}=\widetilde{B}^s_{pq}(L_0)$ is defined as the set of all $f\in \cD'$ such that
\begin{equation}\label{Besov-2}
\|f\|_{\widetilde{B}^s_{pq}}
:= \Big(\sum_{j\ge 0}\big(\|V(\cdot, 2^{-j})^{-s/2}\varphi_j(\sqrt{L_0})f(\cdot)\|_p\big)^q\Big)^{1/q} < \infty.
\end{equation}
Above the usual modification is applied when $p=\infty$ or $q=\infty$.
\end{defn}

\smallskip

We also introduce Besov spaces associated to $L_\cc$:
(a) Classical Besov spaces $\cB^s_{pq}=\cB^s_{pq}(L_\cc)$, and
(b) Nonclassical Besov spaces $\widetilde{\cB}^s_{pq}=\widetilde{\cB}^s_{pq}(L_\cc)$.

\begin{defn}\label{def:Besov-c}
Let $s\in\RR$ and $0< p, q\le \infty$.
Let $\varphi_j$ be as in Definition~\ref{def:Besov-0}.

$(a)$ The Besov space $\cB^s_{pq}=\cB^s_{pq}(L_\cc)$ is defined as the set of all $f\in \cD'$ such that
\begin{equation}\label{Besov-11}
\|f\|_{\cB^s_{pq}}
:= \Big(\sum_{j\ge 0}\big(2^{sj}\|\varphi_j(\sqrt{L_\cc})f\|_p\big)^q\Big)^{1/q} < \infty.
\end{equation}

$(b)$ The Besov space $\widetilde{\cB}^s_{pq}=\widetilde{\cB}^s_{pq}(L_\cc)$ is defined as the set of all $f\in \cD'$ such that
\begin{equation}\label{Besov-22}
\|f\|_{\widetilde{\cB}^s_{pq}}
:= \Big(\sum_{j\ge 0}\big(\|V(\cdot, 2^{-j})^{-s/2}\varphi_j(\sqrt{L_\cc})f(\cdot)\|_p\big)^q\Big)^{1/q} < \infty.
\end{equation}
Above the usual modification is applied when $p=\infty$ or $q=\infty$.
\end{defn}

Several remarks are in order, see \cite{KP}.

(1) The definitions of the spaces
$B^s_{pq}$, $\widetilde{B}^s_{pq}$ or
$\cB^s_{pq}$, $\widetilde{\cB}^s_{pq}$
are independent of the particular selection of the functions $\varphi_0, \varphi$ %or $\Phi_0, \Phi$
with properties \eqref{phi-0}-\eqref{phi}. % or \eqref{phi-00}-\eqref{phi-11}.

(2) The spaces
$B^s_{pq}$, $\widetilde{B}^s_{pq}$ or
$\cB^s_{pq}$, $\widetilde{\cB}^s_{pq}$
are quasi-Banach spaces (Banach spaces if $p, q\ge 1$).

(3) The spaces
$B^s_{pq}$, $\widetilde{B}^s_{pq}$ or
$\cB^s_{pq}$, $\widetilde{\cB}^s_{pq}$
are continuously embedded in $\cD'$.

%We will focus on the spaces $B^s_{pq}$ and $\cB^s_{pq}$.

\subsection{Main result on Besov spaces}\label{subsec:main}

\begin{thm}\label{thm:Besov}
Let $s\in\RR$ and $0< p, q\le \infty$.
Then we have:

$(a)$
$B^s_{pq}(L_0)=\cB^s_{pq}(L_\cc)$ with equivalent norms:
\begin{equation}\label{Besov-equiv}
\|f\|_{B^s_{pq}(L_0)} \sim \|f\|_{\cB^s_{pq}(L_\cc)},\quad f\in B^s_{pq}(L_0).
\end{equation}

$(b)$
$\widetilde{B}^s_{pq}(L_0)=\widetilde{\cB}^s_{pq}(L_\cc)$ with equivalent norms:
\begin{equation}\label{Besov-equiv-2}
\|f\|_{\widetilde{B}^s_{pq}(L_0)} \sim \|f\|_{\widetilde{\cB}^s_{pq}(L_\cc)},\quad f\in B^s_{pq}(L_0).
\end{equation}
\end{thm}

To prove this theorem we need some preparation.

\begin{thm}\label{thm:L0-Lc}
Let $F\in \cS(\RR)$ be even and real-valued.
Then for any $\sigma>0$ and $k\in\NN$ there exists a constant $c(\sigma, k)>0$ such that
\begin{equation}\label{L0-Lc-1}
|L_{0,x}^k \KK_{F(\delta\sqrt{L_\cc})}(x, y)|
\le \frac{c(\sigma, k)\delta^{-2k}}{V(y, \delta)\Big(1+\frac{\rho(x, y)}{\delta}\Big)^{\sigma}},
\quad x, y\in [-1, 1], \; 0<\delta \le 1,
\end{equation}
and
\begin{equation}\label{L0-Lc-2}
|L_{\cc, x}^k \KK_{F(\delta\sqrt{L_0})}(x, y)|
\le \frac{c(\sigma, k)\delta^{-2k}}{V(y, \delta)\Big(1+\frac{\rho(x, y)}{\delta}\Big)^{\sigma}},
\quad x, y\in [-1, 1], \; 0<\delta \le 1.
\end{equation}
Because of the symmetry of the kernels
the above inequalities are valid when $L_{0,x}^k$, $L_{\cc,x}^k$ are replaced by $L_{0,y}^k$, $L_{\cc,y}^k$.
Recall that $L_{0,x}^k \KK(x, y)$ means that the operator $L_{0}^k$ acts on $\KK(x, y)$ as a function of~$x$;
the notation $L_{\cc,x}^k \KK(x, y)$ has a similar meaning.
Also, because $\sigma$ can be arbitrarily large, on account of \eqref{double-3}
the above inequalities are valid in either case with $V(y, \delta)$ replaced by $V(x, \delta)$.
\end{thm}

\begin{proof}
We focus on the proof of \eqref{L0-Lc-1}.
From \eqref{eigenf-L} and \eqref{rep-L0-LcD} it follows that
\begin{align*}
L_{0}^k \psi_n(x)
&= \sum_{\nu=0}^k \gamma_\nu(x) L_\cc^\nu\psi_n(x)
+ \sum_{\nu=0}^{k-1} \eta_\nu(x) \fD L_\cc^\nu\psi_n(x)
\\
&= \sum_{\nu=0}^k \gamma_\nu(x) \chi_n^\nu\psi_n(x)
+ \sum_{\nu=0}^{k-1} \eta_\nu(x) \chi_n^\nu\fD \psi_n(x),
\quad k\in\NN,
\end{align*}
where
$\fD:=(1-x^2)\frac{d}{dx}$,
$\gamma_k(x)=1$,
$\gamma_\nu(x)$ and $\eta_\nu(x)$, $\nu=0, \dots, k-1$, are polynomials of degree $\le 3k$.
In turn, this, \eqref{kernel-fL}, and the rapid decay of $\{F(\delta\sqrt{\chi_n})\}$ imply
%that for any $k\in\NN$
\begin{align*}
L_{0, x}^k \KK_{F(\delta\sqrt{L_\cc})}(x, y)
&= \sum_{\nu=0}^k \gamma_\nu(x) \sum_{n\ge 0} \chi_n^\nu F(\delta\sqrt{\chi_n})\psi_n(x)\psi_n(y)
\\
&+ \sum_{\nu=0}^{k-1} \eta_\nu(x) \fD_x \sum_{n\ge 0} \chi_n^\nu F(\delta\sqrt{\chi_n})\psi_n(x)\psi_n(y)
%\quad k\in\NN,
\end{align*}
and hence
\begin{equation}\label{L0k-f}
L_{0, x}^k \KK_{F(\delta\sqrt{L_\cc})}(x, y)
= \sum_{\nu=0}^k \gamma_\nu(x) \KK_{L_\cc^\nu F(\delta\sqrt{L_\cc})}(x, y)
+ \sum_{\nu=0}^{k-1} \eta_\nu(x) \fD_x \KK_{L_\cc^\nu F(\delta\sqrt{L_\cc})}(x, y),
\end{equation}
where
$\fD_x:=(1-x^2)\partial_x$.

Consider the function
$H(\lambda):= \lambda^{2\nu}F(\lambda)$.
Clearly, $H(\delta\sqrt{L_\cc})= \delta^{2\nu}L_\cc^\nu F(\delta\sqrt{L_\cc})$.
On the other hand, $H\in \cS(\RR)$, $H$ is even and real-valued.
Therefore, by Theorem~\ref{thm:ker-local-S}, applied to $H$, it follows that
for any $\sigma>0$
\begin{equation}\label{local-h}
|\KK_{L_\cc^\nu F(\delta\sqrt{L_\cc})}(x, y)|
\le \frac{c_\sigma\delta^{-2\nu}}{V(y, \delta)\Big(1+\frac{\rho(x, y)}{\delta}\Big)^\sigma}.
\end{equation}
Just as in \eqref{dx-int} we have
\begin{equation}\label{est-DLc-1}
\fD_x \KK_{L_\cc^\nu F(\delta\sqrt{L_\cc})}(x, y)
= \int_x^1 \big[\KK_{L_\cc^{\nu+1} F(\delta\sqrt{L_\cc})}(u, y)
-\cc^2u^2 \KK_{L_\cc^\nu F(\delta\sqrt{L_\cc})}(u, y)\big]du
\end{equation}
and symmetrically
\begin{equation}\label{est-DLc-2}
\fD_x \KK_{L_\cc^\nu F(\delta\sqrt{L_\cc})}(x, y)
= \int_{-1}^x \big[- \KK_{L_\cc^{\nu+1} F(\delta\sqrt{L_\cc})}(u, y)
+ \cc^2u^2 \KK_{L_\cc^\nu F(\delta\sqrt{L_\cc})}(u, y)\big]du.
\end{equation}
Assume $-1\le y < x\le 1$. By \eqref{local-h} and \eqref{est-DLc-1} we obtain
using the change of variables:
$x:=\cos \theta$, $y:=\cos \phi$, $x:=\cos \alpha$,
\begin{align*}
\big|\fD_x &\KK_{L_\cc^\nu F(\delta\sqrt{L_\cc})}(x, y)\big|
\le \int_x^1 \frac{c_\sigma\delta^{-2\nu-2}}{V(y, \delta)\Big(1+\frac{\rho(u, y)}{\delta}\Big)^\sigma} du
\\
&=\frac{c_\sigma\delta^{-2\nu-2}}{V(y, \delta)}
\int_0^\theta \frac{\sin\alpha}{\Big(1+\frac{\phi-\alpha}{\delta}\Big)^\sigma} d\alpha
\le \frac{c_\sigma\delta^{-2\nu-2}}{V(y, \delta)}
\int_0^\theta \frac{1}{\Big(1+\frac{\phi-\alpha}{\delta}\Big)^\sigma} d\alpha.
\end{align*}
Hence
\begin{equation}\label{est-DLc}
\big|\fD_x \KK_{L_\cc^\nu F(\delta\sqrt{L_\cc})}(x, y)\big|
\le \frac{c\delta^{-2\nu-1}}{V(y, \delta)\Big(1+\frac{\rho(x, y)}{\delta}\Big)^{\sigma-1}}.
\end{equation}
In the case when $-1\le x < y\le 1$ we obtain \eqref{est-DLc} as above using \eqref{est-DLc-2}.
Finally, combining \eqref{L0k-f} with \eqref{local-h} and \eqref{est-DLc}
leads to \eqref{L0-Lc-1}.

The proof of \eqref{L0-Lc-2} is similar to the proof of \eqref{L0-Lc-1}.
We omit it.
\end{proof}

We will need the following useful

\begin{lem}\label{lem:integral}
$(a)$ If $\sigma >2$ and $\delta>0$, then
\begin{equation}\label{int-1}
\int_{-1}^1(1+\delta^{-1}\rho(x, y))^{-\sigma}dy \le cV(x, \delta),
\quad x\in [-1, 1].
\end{equation}

$(b)$ If $\sigma >2$, then for $x, y\in[-1, 1]$ and $\delta>0$,
\begin{align}\label{int-2}
\int_{-1}^1 \frac{1}{(1+\delta^{-1}\rho(x, u))^{\sigma}(1+\delta^{-1}\rho(y, u))^{\sigma}}dy
&\le c\frac{V(x, \delta)+V(y, \delta)}{(1+\delta^{-1}\rho(x, y))^{\sigma}}
\\
& \le \frac{cV(x, \delta)}{(1+\delta^{-1}\rho(x, y))^{\sigma-2}}. \nonumber
\end{align}
\end{lem}

For the proof of this lemma, see \cite[Lemma~2.1]{KP}.

\begin{proof}[Proof of Theorem~\ref{thm:Besov}]
We adhere to the notation from Subsection~\ref{subsec:Besov-L0}.

Assume first that $f\in \cB^s_{pq}(L_\cc)$.
We will show that $f\in B^s_{pq}(L_0)$ and
\begin{equation}\label{Besov-est-1}
\|f\|_{B^s_{pq}(L_0)} \le c\|f\|_{\cB^s_{pq}(L_\cc)}.
\end{equation}
Let $\Phi_0, \Phi \in C^\infty(\RR)$ be even, real-valued functions satisfying \eqref{phi-0}-\eqref{phi}.
Just as in Definition~\ref{def:Besov-c}, set $\Phi_j(\lambda):=\Phi(2^{-j}\lambda)$, $j\ge 1$.
It is easy to see (e.g. \cite{FJ2}) that there exist functions $\Psi_0, \Psi\in C^\infty_0(\RR)$
that are even, real-valued, satisfy \eqref{phi-0}-\eqref{phi}, and
\begin{equation}\label{Psi-Phi}
\Psi_0(\lambda)\Phi_0(\lambda) + \sum_{m\ge 1} \Psi(2^{-m}\lambda)\Phi(2^{-m}\lambda) = 1,
\quad \lambda\in \RR.
\end{equation}
Set $\Psi_m(\lambda):=\Psi(2^{-m}\lambda)$, $m\ge 1$.
Then
$\sum_{j\ge 0} \Psi_j(\lambda)\Phi_j(\lambda) = 1$.
From Proposition~\ref{prop:converge} it readily follows that
\begin{equation*}
f= \sum_{m\ge 0} \Psi_m(\sqrt{L_\cc})\Phi_m(\sqrt{L_\cc})f
\quad \hbox{(convergence in $\cD'$).}
\end{equation*}
Let $\varphi_j(\lambda):=\varphi(2^{-j}\lambda)$ be the functions from the definition of
the Besov spaces $B^s_{pq}(L_0)$ in \S \ref{subsec:Besov-L0}.
Then for any $j\ge 0$ we have
\begin{align}\label{rep-phi-j}
\varphi_j(\sqrt{L_0})f &= \sum_{m\ge 0} \varphi_j(\sqrt{L_0})\Psi_m(\sqrt{L_\cc})\Phi_m(\sqrt{L_\cc})f \nonumber
\\
&= \sum_{m\ge 0} \int_I\int_I \KK_{\varphi_j(\sqrt{L_0})}(x, y)
\KK_{\Psi_m(\sqrt{L_\cc})}(y, z)dy \,\Phi_m(\sqrt{L_\cc})f(z)dz,
\end{align}
where $I:=[-1, 1]$.
Two cases present themselves here.

\smallskip

{\em Case 1:} $m\le j$. Assume $j>0$; the case $m=j=0$ is trivial.
Consider the function $H(\lambda):=\varphi(\lambda)\lambda^{-2k}$, $k\in \NN$.
We have
\begin{align*}
\KK_{\varphi_j(\sqrt{L_0})}(x, y)
&=\sum_{n\ge 0}\varphi\big(2^{-j}\sqrt{n(n+1)}\big)\bar{P}_n(x)\bar{P}_n(y)
\\
&=\sum_{n\ge 0}\varphi\big(2^{-j}\sqrt{n(n+1)}\big)[n(n+1)]^k[n(n+1)]^{-k}\bar{P}_n(x)\bar{P}_n(y)
\\
%&= 2^{-2kj}\KK_{L_0^k H(2^{-j}\sqrt{L_0})}(x, y)
&= 2^{-2kj}L_{0, y}^k \KK_{H(2^{-j}\sqrt{L_0})}(x, y).
\end{align*}
Therefore,
\begin{align}\label{case1-1}
\int_I \KK_{\varphi_j(\sqrt{L_0})}(x, y)& \KK_{\Psi_m(\sqrt{L_\cc})}(y, z)dy \nonumber
\\
&= 2^{-2kj}\int_I L_{0, y}^k \KK_{H(2^{-j}\sqrt{L_0})}(x, y) \KK_{\Psi_m(\sqrt{L_\cc})}(y, z)dy
\\
& = 2^{-2kj}\int_I \KK_{H(2^{-j}\sqrt{L_0})(x, y)} L_{0,y}^k \KK_{\Psi_m(\sqrt{L_\cc})}(y, z)dy, \nonumber
\end{align}
where for the last identity we used that the operator $L_0$ is symmetric.
The function $H\in C^\infty(\RR)$ is even, real-valued, and $\supp H\subset [-2, 2]$.
Then by the analogue of Theorem~\ref{thm:ker-local-S} for $L_0$ for any $\sigma>0$
\begin{equation}\label{case1-2}
|\KK_{H(2^{-j}\sqrt{L_0})}(x, y)| \le \frac{c_\sigma}{V(x, 2^{-j})\big(1+2^j\rho(x, y)\big)^{\sigma}}
\le \frac{c2^{2(j-m)}}{V(x, 2^{-m})\big(1+2^m\rho(x, y)\big)^{\sigma}},
\end{equation}
where for the last inequality we used that
$V(x, 2^{-j})\ge c2^{-2(j-m)}V(x, 2^{-m})$, which follows from \eqref{double-2}.
On the other hand, using Theorem~\ref{thm:L0-Lc} we have
\begin{equation}\label{case1-3}
|L_{0,y}^k\KK_{\Psi(2^{-m}\sqrt{L_\cc})}(x, y)|
\le \frac{c(\sigma, k)2^{2km}}{V(x, 2^{-m})\big(1+2^m\rho(x, y)\big)^{\sigma}}.
\end{equation}
Note that the constants $\sigma>0$ and $k\ge 1$ here can be arbitrarily large.
We choose $t$, $k$, and $\sigma$ so that
\begin{equation}\label{sigma-k}
0<t< p,
\quad
2k-2/t-5 > |s|,
\quad
\sigma := 2/t + 5.
\end{equation}
From \eqref{case1-1}-\eqref{case1-3} we infer that
\begin{align*}
\Big|\int_I \KK_{\varphi_j(\sqrt{L_0})}(x, y)& \KK_{\Psi_m(\sqrt{L_\cc})}(y, z)dy\Big|
\\
&\le \frac{c2^{-(2k-2)(j-m)}}{V(x, 2^{-m})^2}\int_I \frac{1}{\big(1+2^m\rho(x, y)\big)^{\sigma}\big(1+2^m\rho(y, z)\big)^{\sigma}} dy
\\
&  \le \frac{c2^{-(2k-2)(j-m)}}{V(x, 2^{-m})\big(1+2^m\rho(x, z)\big)^{\sigma-2}}.
\end{align*}
Here for the last inequality we used \eqref{int-2}.
From this, \eqref{rep-phi-j}, and Lemma~\ref{lem:peetre} it follows that
\begin{align}\label{integ-phi-j-1}
|\varphi_j(\sqrt{L_0})&\Psi_m(\sqrt{L_\cc})\Phi_m(\sqrt{L_\cc})f(x)|
\le \frac{c2^{-(2k-2)(j-m)}}{V(x, 2^{-m})}
\int_I\frac{|\Phi_m(\sqrt{L_\cc})f(z)|}{\big(1+2^m\rho(x, z)\big)^{\sigma-2}} dz \nonumber
\\
&\le c2^{-(2k-2)(j-m)} \sup_{z\in I}\frac{|\Phi_m(\sqrt{L_\cc})f(z)|}{\big(1+2^m\rho(x, z)\big)^{2/t}}
\int_I\frac{V(x, 2^{-m})^{-1}}{\big(1+2^m\rho(x, z)\big)^3} dz
\\
&\le c2^{-(2k-2)(j-m)}\cM_t\big(\Phi_m(\sqrt{L_\cc})f(\cdot)\big)(x), \nonumber
\end{align}
where we used \eqref{int-1} and Lemma~\ref{lem:peetre}.
In turn this and the maximal inequality \eqref{vector-max}, applied to a single function, yields
\begin{align}\label{norm-phi-j1}
\big\|\varphi_j(\sqrt{L_0})\Psi_m(\sqrt{L_\cc})\Phi_m(\sqrt{L_\cc})f\big\|_p
&\le c2^{-(2k-2)(j-m)}\|\cM_t\big(\Phi_m(\sqrt{L_\cc})f(\cdot)\big)\|_p \nonumber
\\
&\le c2^{-(2k-2)(j-m)}\|\Phi_m(\sqrt{L_\cc})f\|_p.
\end{align}

\smallskip

{\em Case 2:} $m>j$.
Let $H(\lambda):= \Psi(\lambda)\lambda^{-k}$.
We write
\begin{align*}
\KK_{\Psi_m(\sqrt{L_\cc})}(y, z)
& =\sum_{n: 1/2<2^{-m}\sqrt{\chi_n}<2} \Psi(2^{-m}\sqrt{\chi_n})\psi_n(y)\psi_n(z)
\\
& =\sum_{n: 1/2<2^{-m}\sqrt{\chi_n}<2} \Psi(2^{-m}\sqrt{\chi_n})\chi_n^k \chi_n^{-k}\psi_n(y)\psi_n(z)
\\
&= 2^{-2km}L_{\cc, y}^k \KK_{H(2^{-m}\sqrt{L_\cc})}(y, z).
\end{align*}
Hence,
\begin{align*}
\int_I \KK_{\varphi_j(\sqrt{L_0})}(x, y)& \KK_{\Psi_m(\sqrt{L_\cc})}(y, z)dy
\\
&= 2^{-2km}\int_I \KK_{\varphi_j(\sqrt{L_0})}(x, y) L_{\cc, y}^k \KK_{H(2^{-m}\sqrt{L_\cc})}(y, z)dy
\\
&= 2^{-2km}\int_I L_{\cc, y}^k \KK_{\varphi_j(\sqrt{L_0})}(x, y) H(2^{-m}\sqrt{L_\cc})(y, z)dy,
\end{align*}
where we used the symmetry of the operator $L_\cc$.
The function $H\in C^\infty(\RR)$ is even, real-valued, and $\supp H\subset [-2, 2]$.
In light of Theorem~\ref{thm:ker-local-S} for any $\sigma>0$
\begin{equation*}
|\KK_{H(2^{-m}\sqrt{L_\cc})}(y, z)| \le \frac{c_\sigma}{V(z, 2^{-m})\big(1+2^m\rho(y, z)\big)^{\sigma}}.
\end{equation*}
On the other hand, from Theorem~\ref{thm:L0-Lc} we have
\begin{equation*}
|L_{\cc, y}^k \KK_{\varphi(2^{-j}\sqrt{L_0})}(x, y)|
\le \frac{c(\sigma, k)2^{2kj}}{V(x, 2^{-j})\big(1+2^j\rho(x, y)\big)^{\sigma}}.
\end{equation*}
Therefore,
\begin{align*}
\Big|\int_I \KK_{\varphi_j(\sqrt{L_0})}(x, y)& \KK_{\Psi_m(\sqrt{L_\cc})}(y, z)dy\Big|
\\
&\le \frac{c2^{-(2k-2)(m-j)}}{V(x, 2^{-j})V(z, 2^{-m})}
\int_I \frac{1}{\big(1+2^j\rho(x, y)\big)^{\sigma}\big(1+2^j\rho(y, z)\big)^{\sigma}} dy
\\
&  \le \frac{c2^{-(2k-2)(m-j)}}{V(z, 2^{-m})\big(1+2^j\rho(x, z)\big)^{\sigma-2}}.
\end{align*}
Here we used again \eqref{int-2}.
We now choose $t$, $k$, and $\sigma$ so that
\begin{equation}\label{sigma-k-2}
0<t<p,
\quad
2k-2/t-5 > |s|,
\quad
\sigma := 2/t+5.
\end{equation}
We use the above, Lemma~\ref{lem:peetre}, and \eqref{int-1} to obtain
\begin{align}\label{integ-phi-j-2}
|\varphi_j(\sqrt{L_0})&\Psi_m(\sqrt{L_\cc})\Phi_m(\sqrt{L_\cc})f(x)|
\le \frac{c2^{-(2k-2)(m-j)}}{V(x, 2^{-m})}
\int_I\frac{|\Phi_m(\sqrt{L_\cc})f(z)|}{\big(1+2^j\rho(x, z)\big)^{\sigma-2}} dz \nonumber
\\
&\le \frac{c2^{-(2k-2)(m-j)}}{V(x, 2^{-m})}
\int_I\frac{2^{(m-j)(2/t+3)}|\Phi_m(\sqrt{L_\cc})f(z)|}{\big(1+2^m\rho(x, z)\big)^{2/t+3}} dz
\\
&\le c2^{-(2k-2/t-5)(m-j)} \sup_{z\in I}\frac{|\Phi_m(\sqrt{L_\cc})f(z)|}{\big(1+2^m\rho(x, z)\big)^{2/t}}
\int_I\frac{V(x, 2^{-m})^{-1}}{\big(1+2^m\rho(x, z)\big)^3} dz \nonumber
\\
&\le c2^{-(2k-2/t-5)(m-j)}\cM_t\big(\Phi_m(\sqrt{L_\cc})f(\cdot)\big)(x). \nonumber
\end{align}
Now, applying the maximal inequality \eqref{vector-max} $(0<t<p)$ for a single function we get
\begin{align}\label{norm-phi-j2}
\big\|\varphi_j(\sqrt{L_0})\Psi_m(\sqrt{L_\cc})\Phi_m(\sqrt{L_\cc})f\big\|_p
&\le c2^{-(2k-2/t-5)(m-j)}\|\cM_t\big(\Phi_m(\sqrt{L_\cc})f(\cdot)\big)\|_p \nonumber
\\
&\le c2^{-(2k-2/t-5)(m-j)}\|\Phi_m(\sqrt{L_\cc})f\|_p.
\end{align}

From \eqref{rep-phi-j}, \eqref{norm-phi-j1}, and \eqref{norm-phi-j2} we obtain
\begin{align}
\|\varphi_j(\sqrt{L_0})f\|_p
&\le c\sum_{m\ge 0}2^{-(2k-2/t-5)|m-j|}\|\Phi_m(\sqrt{L_\cc})f\|_p,
\quad 1\le p\le \infty, \label{norm-phi-f1}
\\
\|\varphi_j(\sqrt{L_0})f\|_p^p
&\le c\sum_{m\ge 0}2^{-(2k-2/t-5)|m-j|p}\|\Phi_m(\sqrt{L_\cc})f\|_p^p,
\quad 0< p< 1. \label{norm-phi-f2}
\end{align}

To complete the proof of \eqref{Besov-est-1} we need the following well known discrete Hardy inequalities:
If $\beta>0$, $0<q<\infty$, and $a_m\ge 0$, then
\begin{align}
\Big(\sum_{j\ge 0}\Big(\sum_{m=0}^j 2^{-\beta(j-m)}a_m\Big)^q\Big)^{1/q}
&\le c\Big(\sum_{m\ge 0}a_m^q\Big)^{1/q}
\quad\hbox{and} \label{hardy-1}
\\
\Big(\sum_{j\ge 0}\Big(\sum_{m\ge j} 2^{-\beta(m-j)}a_m\Big)^q\Big)^{1/q}
&\le c\Big(\sum_{m\ge 0}a_m^q\Big)^{1/q}.  \label{hardy-2}
\end{align}

We next complete the proof of \eqref{Besov-est-1} in the case when $1\le p\le \infty$ and $q<\infty$.
We use \eqref{Besov-1}, \eqref{norm-phi-f1}, \eqref{hardy-1}, \eqref{hardy-2},
and the fact that $2k-2/t-5 > |s|$ to obtain
\begin{align*}
\|f\|_{B^s_{pq}(L_0)}
&= \Big(\sum_{j\ge 0}\big(2^{sj}\|\varphi_j(\sqrt{L_0})f\|_p\big)^q\Big)^{1/q}
\\
&\le c\Big(\sum_{j\ge 0}
\Big(\sum_{m=0}^j2^{-(2k-2/t-5 -s)(j-m)} 2^{sm}\|\Phi_m(\sqrt{L_\cc})f\|_p\Big)^q\Big)^{1/q}
\\
&  + c\Big(\sum_{j\ge 0}
\Big(\sum_{m\ge j+1}2^{-(2k-2/t-5 +s)(m-j)} 2^{sm}\|\Phi_m(\sqrt{L_\cc})f\|_p\Big)^q\Big)^{1/q}
\\
& \le c\Big(\sum_{m\ge 0}\big(2^{sm}\|\Phi_m(\sqrt{L_\cc})f\|_p\big)^q\Big)^{1/q}
\le c\|f\|_{\cB^s_{pq}(L_\cc)}.
\end{align*}
Above we have also used the fact that the definition of $\cB^s_{pq}(L_\cc)$ is independent of the specific selection of
the functions $\varphi_0$, $\varphi$ in Definitions~\ref{def:Besov-c}.

We proceed similarly in the case when $0<p<1$ and $q<\infty$
using \eqref{norm-phi-f2} instead of \eqref{norm-phi-f1}.
The argument in the case $q=\infty$ is easier; we omit it.

The proof of the estimate
\begin{equation}\label{Besov-est-2}
\|f\|_{\cB^s_{pq}(L_\cc)} \le c\|f\|_{B^s_{pq}(L_0)},
\quad f\in B^s_{pq}(L_0),
\end{equation}
is carried out exactly as the proof of \eqref{Besov-est-1},
where the roles of $L_0$, $L_\cc$ and $\varphi_j(\sqrt{L_0})$, $\Phi_m(\sqrt{L_\cc})$
are interchanged.
We omit the details.

The proof of equivalence \eqref{Besov-equiv-2} is similar to the above proof and will be omitted.
\end{proof}

\subsection{Triebel-Lizorkin spaces generated by the operators $L_0$ and $L_\cc$}\label{subsec:F-spaces}

Again following \cite{KP} we next introduce two kinds of Triebel-Lizorkin spaces associated with the operator $L_0$:
(a) Classical Triebel-Lizorkin spaces $F^s_{pq}=F^s_{pq}(L_0)$, and
(b) Nonclassical Triebel-Lizorkin spaces $\widetilde{F}^s_{pq}=\widetilde{F}^s_{pq}(L_0)$.

\begin{defn}\label{def:F-0}
Let $s\in\RR$, $0< p <\infty$, and $0< q\le \infty$.
Also, let $\varphi_0, \varphi \in C^\infty(\RR)$ be even, real-valued functions
satisfying conditions \eqref{phi-0}-\eqref{phi}.

$(a)$ The Triebel-Lizorkin space $F^s_{pq}=F^s_{pq}(L_0)$ is defined as the set of all $f\in \cD'$ such that
\begin{equation}\label{F-1}
\|f\|_{F^s_{pq}} := \Big\|\Big(\sum_{j\ge 0}\big(2^{sj}|\varphi_j(\sqrt{L_0})f(\cdot)|\big)^q\Big)^{1/q}\Big\|_p < \infty.
\end{equation}

$(b)$ The Triebel-Lizorkin space $\widetilde{F}^s_{pq}=\widetilde{F}^s_{pq}(L_0)$ is defined
as the set of all $f\in \cD'$ such that
\begin{equation}\label{F-2}
\|f\|_{\widetilde{F}^s_{pq}}
:= \Big\|\Big(\sum_{j\ge 0}\big(V(\cdot, 2^{-j})^{-s/2}|\varphi_j(\sqrt{L_0})f(\cdot)|\big)^q\Big)^{1/q}\Big\|_p < \infty.
\end{equation}
Above the $\ell^q$-norm is replaced by the $\ell^\infty$-norm when $q=\infty$.
\end{defn}

We also introduce two kinds of Triebel-Lizorkin spaces associated with the operator $L_\cc$:
(a) Classical Triebel-Lizorkin spaces $\cF^s_{pq}=\cF^s_{pq}(L_\cc)$, and
(b) Nonclassical Triebel-Lizorkin spaces $\widetilde{\cF}^s_{pq}=\widetilde{\cF}^s_{pq}(L_\cc)$.

\begin{defn}\label{def:F-c}
Let $s\in\RR$, $0< p <\infty$, and $0< q\le \infty$.
Let $\varphi_0, \varphi \in C^\infty(\RR)$ be even, real-valued functions
satisfying conditions \eqref{phi-0}-\eqref{phi}.

$(a)$ The Triebel-Lizorkin space $\cF^s_{pq}=\cF^s_{pq}(L_\cc)$ is defined as the set of all $f\in \cD'$ such that
\begin{equation}\label{F-11}
\|f\|_{\cF^s_{pq}}
:= \Big\|\Big(\sum_{j\ge 0}\big(2^{sj}|\varphi_j(\sqrt{L_\cc})f(\cdot)|\big)^q\Big)^{1/q}\Big\|_p < \infty.
\end{equation}

$(b)$ The Triebel-Lizorkin space $\widetilde{\cF}^s_{pq}=\widetilde{\cF}^s_{pq}(L_\cc)$ is defined as the set of all $f\in \cD'$ such that
\begin{equation}\label{F-22}
\|f\|_{\widetilde{\cF}^s_{pq}}
:= \Big\|\Big(\sum_{j\ge 0}\big(V(\cdot, 2^{-j})^{-s/2}|\varphi_j(\sqrt{L_\cc})f(\cdot)|\big)^q\Big)^{1/q}\Big\|_p < \infty.
\end{equation}
Above the $\ell^q$-norm is replaced by the $\ell^\infty$-norm when $q=\infty$.
\end{defn}

The Triebel-Lizorkin spaces introduced above have properties similar to the properties
of Besov spaces alluded to in Subsection~\ref{subsec:Besov-L0}, see also \cite{KP}.

\smallskip

The following analogue of Theorem~\ref{thm:Besov} is valid.

\begin{thm}\label{thm:F-spaces}
Let $s\in\RR$, $0< p < \infty$, and $0< q\le \infty$.
Then we have:

$(a)$
$F^s_{pq}(L_0)=\cF^s_{pq}(L_\cc)$ with equivalent norms:
\begin{equation}\label{F-equiv}
\|f\|_{F^s_{pq}(L_0)} \sim \|f\|_{\cF^s_{pq}(L_\cc)},\quad f\in F^s_{pq}(L_0).
\end{equation}

$(b)$
$\widetilde{F}^s_{pq}(L_0)=\widetilde{\cF}^s_{pq}(L_\cc)$ with equivalent norms:
\begin{equation}\label{F-equiv-2}
\|f\|_{\widetilde{F}^s_{pq}(L_0)} \sim \|f\|_{\widetilde{\cF}^s_{pq}(L_\cc)},\quad f\in \widetilde{F}^s_{pq}(L_0).
\end{equation}
\end{thm}

\begin{proof}
This proof uses the ingredients from the proof of Theorem~\ref{thm:Besov}.

Assuming that $f\in \cF^s_{pq}(L_\cc)$ we will show that $f\in F^s_{pq}(L_0)$ and
\begin{equation}\label{F-est-1}
\|f\|_{F^s_{pq}(L_0)} \le c\|f\|_{\cF^s_{pq}(L_\cc)}.
\end{equation}
Let $\Phi_0, \Phi \in C^\infty(\RR)$ be even, real-valued functions satisfying \eqref{phi-0}-\eqref{phi}.
Set $\Phi_j(\lambda):=\Phi(2^{-j}\lambda)$, $j\ge 1$.
Then just as in the proof of Theorem~\ref{thm:Besov} there exists
even, real-valued functions $\Psi_0, \Psi\in C_0^\infty(\RR)$
that satisfy \eqref{phi-0}, \eqref{phi}, and \eqref{Psi-Phi}.
Set $\Psi_m(\lambda):=\Psi(2^{-m}\lambda)$, $m\ge 1$.
Also, let $\varphi_j(\lambda):=\varphi(2^{-j}\lambda)$ be the functions from Definition~\ref{def:F-0}.

We now proceed as in the proof of Theorem~\ref{thm:Besov} to obtain
\eqref{integ-phi-j-1} and \eqref{integ-phi-j-2}.
From \eqref{integ-phi-j-1} and \eqref{integ-phi-j-2} we infer that
%From these estimates we infer that
\begin{equation*}
|\varphi_j(\sqrt{L_0})\Psi_m(\sqrt{L_\cc})\Phi_m(\sqrt{L_\cc})f(x)|
\le c2^{-(2k-2/t-5)|m-j|}\cM_t\big(\Phi_m(\sqrt{L_\cc})f(\cdot)\big)(x)
\end{equation*}
and hence, in light of \eqref{rep-phi-j},
\begin{equation*}
|\varphi_j(\sqrt{L_0})f(x)|
\le c\sum_{m\ge 0} 2^{-(2k-2/t-5)|m-j|}\cM_t\big(\Phi_m(\sqrt{L_\cc})f(\cdot)\big)(x).
\end{equation*}
Similarly as in the proof of Theorem~\ref{thm:Besov},
we choose $t$ and $k$ so that $0<t< \min\{p, q\}$ and $2k-2/t-5 > |s|$.
We use the above in the definition of $\|f\|_{F^s_{pq}(L_0)}$ to obtain
\begin{align*}
\|f\|_{F^s_{pq}(L_0)} &= \Big\|\Big(\sum_{j\ge 0}\big(2^{sj}|\varphi_j(\sqrt{L_0})f(\cdot)|\big)^q\Big)^{1/q}\Big\|_p
\\
&\le c\Big\|\Big(\sum_{j\ge 0}
\Big[2^{sj} \sum_{m\ge 0} 2^{-(2k-2/t-5)|m-j|}\cM_t\big(\Phi_m(\sqrt{L_\cc})f\big)(\cdot)\Big]^q\Big)^{1/q}\Big\|_p
\\
&\le c\Big\|\Big(\sum_{m\ge 0}
\Big[2^{sm}\cM_t\big(\Phi_m(\sqrt{L_\cc})f\big)(\cdot)\Big]^q\Big)^{1/q}\Big\|_p
\\
&\le c\Big\|\Big(\sum_{m\ge 0}
\big[2^{sm} |\Phi_m(\sqrt{L_\cc})f(\cdot)|\big]^q\Big)^{1/q}\Big\|_p
\le c\|f\|_{\cF^s_{pq}(L_\cc)},
\end{align*}
which confirms \eqref{F-est-1}.
Above for the second inequality we used the Hardy inequalities \eqref{hardy-1}-\eqref{hardy-2}
and the fact that $2k-2/t-5 >|s|$,
for the third inequality we used the maximal inequality \eqref{vector-max} and the fact that $0<t< \min\{p, q\}$,
and for the last inequality we used the fact that the definition of $\cF^s_{pq}$ is independent of the specific selection of
the functions $\varphi_0$, $\varphi$ in Definitions~\ref{def:F-c}.

To prove that
$\|f\|_{\cF^s_{pq}(L_\cc)} \le c\|f\|_{F^s_{pq}(L_0)}$
for $f\in F^s_{pq}(L_0)$,
one simply switches the roles of $L_0$, $L_\cc$ and $\varphi_j(\sqrt{L_0})$, $\Phi_m(\sqrt{L_\cc})$
in the above proof. We omit the details.

The proof of equivalence \eqref{F-equiv-2} is similar to the above proof and will be omitted.
\end{proof}

\section{Frame characterization of Besov and Triebel-Lizorkin spaces}\label{sec:fremes}

The discrete frame characterization of the PSWF Besov and Triebel-Lizorkin spaces is an important part of our theory.
There are two approaches available for construction of frames based on the PSWFs.

One option is to apply the general method developed in \cite{KP}.
All ingredients for this constriction are in place:
a sampling theorem and a quadrature formula can be developed just as in \cite{CKP, KP},
well localized kernels are readily available,
and lower and upper bound estimates on the norms of the kernels can be established easily.
Then the construction of a frame $\{\theta_{\xi}\}_{\xi\in\cX}$ and a dual frame $\{\tilde\theta_{\xi}\}_{\xi\in\cX}$
is straightforward.
As a consequence of the general theory from \cite{KP} these frames will allow to completely characterize
the PSWF Besov and Triebel-Lizorkin spaces.
We skip the further details.

A frame based on Legendre polynomials and
related highly localized kernels like the ones in \S\ref{subsec:smooth-f-c}
has been constructed in \cite{PX}.
Just as in \cite{KP} this frame can be used for characterization of the Besov and Triebel-Lizorkin spaces associated
with the operator $L_0$ and hence with $L_\cc$ from \S\ref{sec:B-F-spaces}.
In turn this polynomial frame can be used as a backbone for the construction of a pair of dual frames
$\{\theta_{\xi}\}_{\xi\in\cX}$, $\{\tilde\theta_{\xi}\}_{\xi\in\cX}$ in terms of the PSWFs
by application of the small perturbation method from \cite{DKKP}.
The highly localized kernels needed for this construction are available from \S \ref{sec:func-calc}.
This pair of dual frames will allow for complete characterization of
the PSWF Besov and Triebel-Lizorkin spaces.
We omit the details.

The construction of a frame based on the PSWFs would be particularly simple and elegant
(like in the polynomial case, see \cite{PX})
if there was a quadrature formula on $[-1, 1]$ with positive weights of the right size
that is exact for all products $g\cdot h$ of functions $g, h\in \Sigma_N^\cc$ (see \eqref{def-Sigma-N}).
The construction of such a quadrature formula is an open problem.

\section{Heat kernel generated by perturbation of the Jacobi operator}\label{sec:Jacobi}

In this and next sections we show that our method for establishing Gaussian bounds for heat kernels
is not limited to the PSWFs of order zero.
We first apply it to perturbations of the Jacobi operator $\LLL$:
\begin{equation}\label{Jacobi-oper}
\LLL f(x):=-\frac{1}{\omega(x)}\frac{d}{dx}\Big[\omega(x)a(x)\frac{df}{dx}(x)\Big],
\quad D(\LLL) := C^2[-1, 1],
\end{equation}
where
\begin{equation*}
\omega(x)=\omega_{\alpha, \beta}(x):= (1-x)^\alpha(1+x)^\beta,
\quad
\alpha, \beta >-1,
\quad a(x):=(1-x^2).
\end{equation*}

Consider now a perturbation of the operator $\LLL$ of the form
\begin{equation}\label{def-L-V}
\LLV f(x):= \LLL f(x)+ V(x)f(x),
\quad \quad D(\LLV) := C^2[-1, 1],
\end{equation}
where $V\in L^\infty[-1, 1]$ and $V\ge 0$. % and $V\not\equiv 0$.

%{\bf (How about assuming that $V\in C^2[-1, 1]$ or $V\in L^\infty[-1, 1]$?)}

From the Sturm-Liouville theory it follows that
$\LLV$ is an essentially self-adjoint positive operator in $L^2(\omeg):=L^2([-1, 1], \omeg(x)dx)$.
Denote by
\begin{equation}\label{eigen-v-V}
\psi_1, \psi_2, \dots
\quad\hbox{and}\quad
0<\chi_1<\chi_2 < \cdots
\end{equation}
the eigenfunctions and eigenvalues of the operator $\LLV$.
We assume that the eigenfunctions $\{\psi_n\}$ are normalized so that
\begin{equation*}
\langle \psi_n, \psi_k \rangle := \int_{-1}^1 \psi_n(x)\psi_k(x)\omeg(x)dx = \delta_{nk}.
\end{equation*}
Observe that $\{\psi_n\}_{n=1}^\infty$ is an orthonormal basis for $L^2(\omeg)$. %$L^2([-1, 1], \omeg(x)dx)$.

We are interested in the space localization of the heat kernel
\begin{equation}\label{Jacobi-Heat}
\KK_{\exp(-t\LLV)}(x, y)= \sum_{n=1}^\infty e^{-\chi_n t}\psi_n(x)\psi_n(y).
\end{equation}

Clearly, the operator $\LLV$ is a perturbation of the Jacobi operator $\LLL$ and
hence $\{\psi_n\}$ are closely related to the Jacobi polynomials.
We denote by $\tP_n$ the $n$th degree Jacobi polynomial with the normalization
\begin{equation*}
\langle \tP_n, \tP_k \rangle := \int_{-1}^1 \tP_n(x)\tP_k(x)\omeg(x)dx = \delta_{nk}.
\end{equation*}
As is well known \cite{Sz} the Jacobi polynomials are eigenfunctions of the Jacobi operator $\LLL$,
i.e.
\begin{equation}\label{LLL-P}
\LLL\tP_n= \lambda_n \tP_n,
\quad
\lambda_n= n(n+\alpha+\beta+1).
\end{equation}

The next theorem relates $\{\chi_n\}$ to $\{\lambda_n\}$.

\begin{thm}\label{thm:eigenvalues}
The eigenvalues of $\LLV$ and $\LLL$ are related by the inequalities
\begin{equation}\label{eigenvalues}
\lambda_n \le \chi_n \le \lambda_n+\|V\|_\infty, \quad n=1, 2, \dots.
\end{equation}
\end{thm}

\begin{proof}
Denote by $\overline{\LLL}$ and $D(\overline{\LLL})$ the closure of
the Jacobi operator $\LLL$ and its domain in $L^2(\omeg)$.
Likewise we let $\overline{\LLV}$ and $D(\overline{\LLV})$ be the closure of $\LLV$ and its domain.
From $\LLV=\LLL+V$ it follows that $D(\overline{\LLV})=D(\overline{\LLL})$
and $\overline{\LLV}=\overline{\LLL}+V$.
Observe that $\overline{\LLL}$ and $\overline{\LLV}$ are positive self-adjoint operators on $L^2(\omeg)$.

By the Poincar\'{e} min-max characterization of the eigenvalues of a positive self-adjoint operator,
applied to $\overline{\LLV}$ first and then to $\overline{\LLL}$ it follows that
\begin{align*}
\chi_n &= \min_{H\subset D(\overline{\LLV}), \, \dim H=n} \max_{f\in H, \|f\|_2=1}
\langle f, \overline{\LLV} f\rangle^2
\\
&\le  \min_{H\subset D(\overline{\LLL}), \, \dim H=n} \max_{f\in H, \|f\|_2=1}
\langle f, \overline{\LLL} f\rangle^2 + \|V\|_\infty\|f\|_2^2
\\
&= \lambda_n+ \|V\|_\infty.
\end{align*}
Above in both cases the minimum is over all $n$-dimensional subspaces $H$.

For the proof of the lower bound estimate in \eqref{eigenvalues}
we observe that from $\overline{\LLV}=\overline{\LLL}+V$ and $V\ge 0$
it follows that
$\langle f, \overline{\LLV} f\rangle^2 \ge \langle f, \overline{\LLL} f\rangle^2$
for all $f\in L^2(\omeg)$,
which implies $\chi_n\ge \lambda_n$
by using the min-max theorem similarly as above.

Observe that in our case of positive self-adjoint operators with discrete spectrums
the min-max theorem is quite obvious;
we will not elaborate any further.
\end{proof}

We consider $[-1, 1]$ equipped with the measure
\begin{equation*}
d\mu(x) := \omeg(x)dx = \omega_{\alpha, \beta}(x)dx
\end{equation*}
and the distance
\begin{equation}\label{distance-J}
\rho(x, y):= |\arccos x - \arccos y|.
\end{equation}
We denote the balls on $[-1, 1]$ by
$B(x, r):=\{y\in [-1, 1]: \rho(x, y)<r\}$
and set
$V(x, r):=\mu(B(x, r))$.
As is shown in \cite[see (7.1)]{CKP}
\begin{equation}\label{size-V}
V(x, r)\sim r(1-x+r^2)^{\alpha+1/2}(1+x+r^2)^{\beta+1/2},
\quad x\in [-1, 1], \; 0<r\le \pi,
\end{equation}
and obviously $V(x, r)=V(x, \pi)\sim 1$ if $r>\pi$.
Hence, the measure $\mu$ has the doubling property, see \eqref{doubling-X}.
%\begin{equation}\label{doubl-J}
%V(x, 2t) \le cV(x, t),
%\quad x\in [-1, 1], t>0.
%\end{equation}

The heat kernel associated to the Jacobi operator $\LLL$ takes the form
\begin{equation}\label{Jacobi-HK}
\KK_{\exp(-t\LLL)}(x, y) = \sum_{n=0}^\infty e^{-\lambda_n t}\tP_n(x)\tP_n(y).
\end{equation}
Gaussian upper and lower bounds for $\KK_{\exp(-t\LLL)}(x, y)$ just as in \eqref{Gauss-L0}
are established in \cite[Theorem 7.2]{CKP}.

%There exist constants $c_1, c_2, c_3, c_4 >0$ such that for all $x,y \in [-1, 1]$ and $t>0$
%\begin{equation}\label{Gauss-Jacobi}
%\frac{c_1\exp\big(- \frac{\rho(x,y)^2}{c_2t}\big)}{\big[V(x, \sqrt t)V(y, \sqrt t)\big]^{1/2}}
%\le  \KK_{\exp(-t\LLL)}(x,y)
%\le \frac{c_3\exp\big(- \frac{\rho(x,y)^2}{c_4t}\big)}{\big[V(x, \sqrt t)V(y, \sqrt t)\big]^{1/2}}.
%\end{equation}

We next establish Gaussian bounds for the heat kernel $\KK_{\exp(-t\LLV)}(x, y)$:

\begin{thm}\label{thm:Gauss-L1}
There exist constants $c_1, c_2, c_3, c_4 >0$ such that for all $x,y \in [-1, 1]$ and $t>0$
\begin{equation}\label{Gauss-L1}
\frac{c_1e^{-t\|V\|_\infty}\exp\big(- \frac{\rho(x,y)^2}{c_2t}\big)}{\big[V(x, \sqrt t)V(y, \sqrt t)\big]^{1/2}}
\le  \KK_{\exp(-t\LLV)}(x,y)
\le \frac{c_3\exp\big(- \frac{\rho(x,y)^2}{c_4t}\big)}{\big[V(x, \sqrt t)V(y, \sqrt t)\big]^{1/2}}.
\end{equation}
\end{thm}

The proof of this theorem relies on our general result from Theorem~\ref{thm:abstract}
and the Gaussian bounds for the heat kernel $\KK_{\exp(-t\LLL)}(x, y)$.
It is carried out just as the proof of Theorem~\ref{thm:main}; we omit it.

\subsection*{H\"{o}lder continuity of the heat kernel \boldmath $\KK_{\exp(-t\LLV)}(x, y)$}

In \cite[Section 7]{CKP} the quadratic form
$
\cE(f, g) := \langle \LLL f, g\rangle
$
is considered and shown that it generates a regular strictly local Dirichlet space.
Furthermore, the local scale-invariant Poincar\'{e} inequality is established,
which along with the doubling property of the measure implies Gaussian bounds %\eqref{Gauss-Jacobi}
for the kernel $\KK_{\exp(-t\LLL)}(x, y)$.
As was explained in \S \ref{subsec:Holder} the Gaussian bounds for $\KK_{\exp(-t\LLL)}(x, y)$ %\eqref{Gauss-Jacobi}
imply the H\"{o}lder continuity of the kernel $\KK_{\exp(-t\LLL)}(x, y)$ in the form of \eqref{Le-Holder}.
%There exists a constant $0<\gamma\le 1$ such that
%\begin{equation}\label{KL-Holder}
%|\KK_{\exp(-t\LLL)}(x, y)-\KK_{\exp(-t\LLL)}(x', y)|
%\le c_5\Big(\frac{\rho(x, x')}{\sqrt{t}}\Big)^\gamma
%\frac{\exp\big(- \frac{\rho(x,y)^2}{c_6t}\big)}{\big[V(x, \sqrt t)V(y, \sqrt t)\big]^{1/2}},
%\end{equation}
%for all $x, x', y \in [-1, 1]$ and $t>0$, whenever $\rho(x, x')\le \sqrt{t}$.
%
In turn, applying Theorem~\ref{thm:Holder}, the Gaussian bounds for $\KK_{\exp(-t\LLL)}(x, y)$
and its H\"{o}lder continuity imply
%estimates \eqref{Gauss-Jacobi} and \eqref{KL-Holder} imply
that $\KK_{\exp(-t\LLV)}(x, y)$ is H\"{o}lder continuous as well:
\begin{align}\label{KLV-Holder}
&|\KK_{\exp(-t\LLV)}(x, y)-\KK_{\exp(-t\LLV)}(x', y)|
\\
&\quad\le c_7(1+t\|V\|_\infty)\Big(\frac{\rho(x, x')}{\sqrt{t}}\Big)^\gamma
\frac{\exp\big(- \frac{\rho(x,y)^2}{c_8t}\big)}{\big[V(x, \sqrt t)V(y, \sqrt t)\big]^{1/2}},
\;\; \hbox{if $\rho(x, x')\le \sqrt{t}$.}\nonumber
\end{align}
%if $\rho(x, x')\le \sqrt{t}$.

\begin{rem}\label{rem:Holder-LV}
A straightforward adaptation of the proof of Theorem~\ref{thm:Holder-PSWF}
shows that estimate \eqref{KLV-Holder} is valid with $\gamma =1$
in the case when $\alpha, \beta\ge -1/2$.
\end{rem}

\subsection*{Generalized prolate spheroidal wave functions}%\label{subsec:general-PSWFs}

The generalized prolate spheroidal wave functions (GPSWFs) in dimension one are introduced in \cite{Wang_Zhang}
as eigenfunctions of the operator
\begin{equation}\label{Gegenbauer-oper}
L_\cc^\alpha :=-\frac{1}{(1-x^2)^\alpha}\frac{d}{dx}\Big[(1-x^2)^{\alpha+1}\frac{d}{dx}\Big] +\cc^2x^2,
\quad \alpha>-1,
\end{equation}
with domain $D(L_\cc^\alpha) := C^2[-1, 1]$,
where $\cc>0$ is a fixed constant.

It was shown in \cite{Wang_Zhang} that the operator $L_\cc^\alpha$ commute with the operator
\begin{equation*}
\cF_\cc^\alpha f(x):=\int_{-1}^1 e^{i\cc xt}f(t) \omeg(t)dt.
\end{equation*}
Denote by $\psi_{1}^\alpha, \psi_{2}^\alpha, \dots$ the joint eigenfunctions of
the operators $L_\cc^\alpha$ and $\cF_\cc^\alpha$.
Clearly, $\{\psi_n^\alpha\}$ is a natural generalization of
the PSWFs $\{\psi_n\}$ (when $\alpha =0$), see the introduction.

From Theorem~\ref{thm:Gauss-L1} it follows that we have Gaussian bounds for the (heat) kernel
$\KK_{\exp(-tL_\cc^\alpha)}(x,y)$ of the semigroup $e^{-tL_\cc^\alpha}$
generated by the operator $L_\cc^\alpha$:
For any $x, y\in [-1, 1]$ and $t>0$
\begin{equation}\label{Gauss-Lc}
\frac{c_1e^{-\cc^2 t}\exp\big(- \frac{\rho(x,y)^2}{c_2t}\big)}{\big[V(x, \sqrt t)V(y, \sqrt t)\big]^{1/2}}
\le  \KK_{\exp(-tL_\cc^\alpha)}(x,y)
\le \frac{c_3\exp\big(- \frac{\rho(x,y)^2}{c_4t}\big)}{\big[V(x, \sqrt t)V(y, \sqrt t)\big]^{1/2}}.
%\quad \forall x,y \in [-1, 1], \; \forall t>0.
\end{equation}

Note that in the case when $\alpha=0$ the above result is
simply Theorem~\ref{thm:main}.

The H\"{o}lder continuity of the heat kernel $\KK_{\exp(-tL_\cc^\alpha)}(x,y)$
follows from \eqref{KLV-Holder}.

\subsection*{Besov and Triebel-Lizorkin spaces associated to the GPSWFs on \boldmath $[-1,1]$}

It should be pointed out that in the setting on $[-1,1]$ with the Jacobi weight,
the distance $\rho(x,y)$ from \eqref{distance-J},
and the perturbation $\LLV$ of the Jacobi operator $\LLL$ or
in the setting of the Generalized prolate spheroidal wave functions
the Gaussian bounds for the heat kernel and its H\"{o}lder continuity
allow for the development of the associated Besov and Triebel-lizorkin spaces
in complete generality (with a complete set of indices).
For this one has to simply apply the general scheme developed in \cite{CKP, KP}.

\section{Prolate spheroidal wave functions on the ball}\label{sec:PSWF-ball}

The prolate spheroidal wave functions on the unit ball $\BB^d$ in $\RR^d$, $d\ge 2$,
were introduced by D. Slepian in \cite{Slep-1}
as a natural generalization of the one-dimensional prolate spheroidal wave functions of order zero.
Slepian's results from \cite{Slep-1} were somewhat further refined in \cite{ZLWZ}.

The PSWFs on $\BB^d$
are defined as the eigenfunctions of the following two commuting operators:
\begin{equation}\label{def-L-ball-intro}
L_\cc f(x):=-\sum_{i=1}^{d}\partial_i^2f(x) + \sum_{i=1}^d \sum_{j=1}^d x_ix_j  \partial_{i}\partial_jf(x)
+ (n+2\gamma)\sum_{i=1}^{d} x_i\partial_if(x) + c^2\|x\|^2f(x)
\end{equation}
and
\begin{equation}\label{def-Fc-B}
F_\cc f(x) :=\int_{\BB^d} f(u)e^{-i\cc x\cdot u}\omega_\gamma(u)du.
\end{equation}
Here $\omega_\gamma(x):= (1-\|x\|^2)^{\gamma-1/2}$, $\gamma >-1/2$, is a weight function,
and $\cc>0$ is a fixed constant.
We use classical notation for the vectors
$x=(x_1, \dots, x_d)\in \RR^d$,
the inner product
$x\cdot y:=\sum_{j=1}^{d}x_jy_j$,
and the Euclidean norm
$\|x\|:=\sqrt{x\cdot x}$.
We consider the unit ball $\BB^d$ equipped with the weighted measure
\begin{equation}\label{mes-ball}
d\mu(x) := \omega_\gamma(x)dx:=(1-\| x\|^2)^{\gamma-1/2} dx,  \quad  \gamma >-1/2,
\end{equation}
and introduce the shorthand notation
$L^2(\omega_\gamma):=L^2(\BB^d, \omega_\gamma(x)dx)$.

Let $\cH_n$ be the set of all spherical harmonics of degree $n$
on the unit sphere $\SS^{d-1}$ in $\RR^d$.
As is well known the dimension of $\cH_n$ is
$N(n, d)= \frac{2n+d-2}{n}\binom{n+d-3}{n-1} \sim n^{d-2}$.
In what follows we will denote by $\{Y_{nj}: 1\le j \le N(n, d)\}$
an orthonormal basis for $\cH_n$ consisting of real-valued functions.

The eigenfunctions of the operators $L_\cc$ take the form
\begin{equation}\label{psi-knj}
\psi_{kj}^n(x)= \|x\|^n\phi_{nk}(2\|x\|^2-1)Y_{nj}\Big(\frac{x}{\|x\|}\Big),
\quad 1\le j\le N(n, d),
\end{equation}
where $\{\phi_{nk}\}$ are the eigenfunctions of a certain differential operator.
For the details, see \cite{ZLWZ}.
Thus,
\begin{equation}\label{PSWF-b}
L_\cc \psi_{kj}^n(x) = \chi_{nk} \psi_{kj}^n(x),\quad n, k\in \NN_0,\; 1\le j\le N(n, d),
\end{equation}
where $\{\chi_{nk}\}$ are the eigenvalues of $L_\cc$.
Note that the PSWFs $\{\psi_{kj}^n\}$ form an orthonormal basis for $L^2(\omega_\gamma)$.

The operator $L_\cc$ is positive and self-adjoint on $L^2(\omega_\gamma)$.
Hence, $L_\cc$ generates a semigroup $\exp(-tL_\cc)$ with
a (heat) kernel $p_t(x, y)=\KK_{\exp(-tL_\cc)}(x, y)$ of the form
\begin{align}\label{heat-ker}
p_t(x, y) &= \sum_{m=0}^\infty \sum_{n+2k=m} e^{-\chi_{nk}t}\sum_{j=1}^{N(n, d)}\psi_{kj}^n(x)\psi_{kj}^n(y)
\\
&= \sum_{n=0}^\infty \sum_{k=0}^\infty e^{-\chi_{nk}t}\sum_{j=1}^{N(n, d)}\psi_{kj}^n(x)\psi_{kj}^n(y). \nonumber
\end{align}

As will be seen below it is natural to use the distance
\begin{equation}\label{dist-ball}
\rho(x,y) := \arccos \big(x\cdot y + \sqrt{1-\| x\|^2}\sqrt{1-\| y\|^2}\big)
\end{equation}
in the current setting.
We will use standard notation for balls $B(x, r)\subset \BB^d$:
\begin{equation}\label{B-ball}
B(x, r):=\{y\in\BB^d: \rho(x, y)<r\} \quad\hbox{and set}\quad V(x, r):=\mu(B(x, r)).
\end{equation}
As is well known (see e.g. \cite[Lemma 11.3.6]{DaiXu})
\begin{equation}\label{V-B}
V(x, r)\sim r^d(1-\|x\|^2+r^2)^\gamma, \quad 0<r\le \pi,
\end{equation}
and obviously $V(x, r)=V(x, \pi) \sim 1$ if $r>\pi$.
Therefore, $\mu$ is a doubling measure.
%\begin{equation}\label{doubling-ball}
%V(x, 2r)\le cV(x, r), \quad x\in\BB^d, \; r>0.
%\end{equation}

\subsection*{Gaussian bounds for the heat kernel associated to $L_\cc$ on $\BB^d$}

\begin{thm}\label{thm:main-ball}
The PSWF heat kernel $p_t(x,y)$ on $\BB^d$ $($see \eqref{heat-ker}$)$ has the following Gaussian upper and lower bounds:
There exist constants $c_1, c_2, c_3, c_4 >0$ such that for all $x,y \in \BB^d$ and $t>0$
\begin{equation}\label{gauss-ball-Lc}
\frac{c_1e^{-t\cc^2}\exp\big(- \frac{\rho(x,y)^2}{c_2t}\big)}{\big[V(x, \sqrt t)V(y, \sqrt t)\big]^{1/2}}
\le  p_t(x,y)
\le \frac{c_3\exp\big(- \frac{\rho(x,y)^2}{c_4t}\big)}{\big[V(x, \sqrt t)V(y, \sqrt t)\big]^{1/2}}.
\end{equation}
\end{thm}

A remark similar to Remark~\ref{rem:main} is valid here.

The proof of Theorem~\ref{thm:main-ball} will again depend on our general result from Theorem~\ref{thm:abstract}.
We will utilize the close relationship of the PSWFs on $\BB^d$ with orthogonal polynomials on $\BB^d$.

\subsection{Orthogonal polynomials on $\BB^d$ and the associated heat kernel}\label{sec:plyn-ball}

Here we first recall the Gaussian localization of the heat kernel
associated to the operator $L_0$ (see \eqref{def-L-ball-intro} with $\cc=0$) on $\BB^d$ from \cite{KPX1, KPX2}.
Note that the operator $L_0$ in the setting described above
is essentially self-adjoint and positive.

Denote by $\cP_m$ the set of all algebraic polynomials of degree $\le m$ in $d$ variables,
and let $\cV_m$ be the orthogonal compliment of $\cP_{m-1}$ to
$\cP_m$ in $L^2(\omega_\gamma)$, $m\ge 1$.
Then $\cP_m=\cP_{m-1}\bigoplus \cV_m$.
Denote $\cV_0:=\cP_0$.
As is well known (see e.g. \cite[\S2.3.2]{DX}) $\cV_m$, $m=0, 1, \dots$,
are eigenspaces of the operator $L_0$,
more precisely,
\begin{equation}\label{eigen-sp}
L_0 P=\lambda_m P,\quad \forall P\in\cV_m,
\;\;\hbox{where $\lambda_m:= m(m+d+2\gamma-1)$, $\; m=0, 1, \dots $}.
\end{equation}
Let $\tilde P_{mj}$, $j=1, \dots, \dim \cV_m$,
be an arbitrary real orthonormal basis for $\cV_m$ in $L^2(\omega_\gamma)$.
Denote $N_m:= \dim \cV_m = \binom{m+d-1}{m}$.
Then
\begin{equation}\label{orth-proj}
P_m(x, y):=\sum_{j=1}^{N_m}\tilde P_{mj}(x)\tilde P_{mj}(y), \quad x,y\in \BB^d,
\end{equation}
is the kernel of the orthogonal projector onto $\cV_m$.
The heat kernel $\KK_{\exp(-tL_0)}(x,y)$, $t>0$, of the semigroup $e^{-tL_0}$ takes the form
\begin{equation}\label{ball-HK}
\KK_{\exp(-tL_0)}(x,y)=\sum_{m=0}^\infty e^{-\lambda_m t} P_m(x, y).
\end{equation}

One of the main results in \cite{KPX1} (see also \cite{KPX2}) asserts that
the heat kernel $\KK_{\exp(-tL_0)}(x, y)$ has Gaussian upper and lower bounds, just as in \eqref{Gauss-L0}.
%that is,
%there exist constants $c_1, c_2, c_3, c_4 >0$ such that for all $x,y \in \BB^d$ and $t>0$
%\begin{equation}\label{gauss-ball}
%\frac{c_1\exp\big(- \frac{\rho(x,y)^2}{c_2t}\big)}{\big[V(x, \sqrt t)V(y, \sqrt t)\big]^{1/2}}
%\le  \KK_{\exp(-tL_0)}(x,y)
%\le \frac{c_3\exp\big(- \frac{\rho(x,y)^2}{c_4t}\big)}{\big[V(x, \sqrt t)V(y, \sqrt t)\big]^{1/2}}.
%\end{equation}

\subsection*{Orthogonal polynomial basis on the ball: Zernike polynomials}

One specific polynomial basis on the ball (Zernike polynomials)
plays an important role in developing the PSWFs on the ball.

We will denote by $\PP_n^{(\alpha, \beta)}(u)$ the $n$th degree Jacobi polynomial
and will assume that the Jacobi polynomials are normalized so that
\begin{equation}\label{jacobi-norm}
\int_{-1}^1 \PP_n^{(\alpha, \beta)}(u)\PP_k^{(\alpha, \beta)}(u)\omega_{\alpha, \beta}(u) du = 2^{\alpha+\beta+2}\delta_{nk},
\end{equation}
where $\omega_{\alpha, \beta}(u):=(1-u)^\alpha(1+u)^\beta$, $\alpha, \beta>-1$.

As before we assume that $\{Y_{nj}: 1\le j \le N(n, d)\}$
is an orthonormal basis for the set $\cH_n$ of spherical harmonics of degree $n$ consisting of real-valued functions.

One particular orthonormal basis for $\cV_m$ is given by the Zernike polynomials (see e.g. \cite{DaiXu, DX}):
\begin{equation*}%\label{pol-basis}
P_{kj}^n(x) := \PP_k^{(\gamma-\frac{1}{2}, n+\frac{d}{2}-1)}(2\|x\|^2-1)\|x\|^n Y_{nj}\Big(\frac{x}{\|x\|}\Big),
\;\; n+2k=m,\; 1\le j\le N(n, d).
\end{equation*}

In terms of the polynomial basis $\{P_{kj}^n\}$ the heat kernel associated to $L_0$ takes the form:
\begin{equation}\label{heat-ker-P}
\KK_{\exp(-tL_0)}(x,y) = \sum_{m=0}^\infty e^{-\lambda_m t}\sum_{n+2k=m}\sum_{j=1}^{N(n, d)} P_{kj}^n(x)P_{kj}^n(y).
\end{equation}

The close relationship between the PSWFs $\{\psi_{kj}^n\}$ and the orthogonal polynomials $\{P_{kj}^n\}$ on $\BB^d$
is reflected in the following relation between the respective eigenvalues (see \cite[Theorem 3.2]{ZLWZ}):
\begin{equation}\label{chi-lam-1}
\lambda_{n+2k} < \chi_{nk} < \lambda_{n+2k} + \cc^2,
\end{equation}
implying
\begin{equation}\label{chi-lam-2}
(n+2k)(n+2k+d+2\gamma-1) < \chi_{nk} < (n+2k)(n+2k+d+2\gamma-1) +\cc^2.
\end{equation}

\subsection{Proof of Theorem~\ref{thm:main-ball}}\label{subsec:proof-main}

As was indicated above the proof of Theorem~\ref{thm:main-ball} will depend on the general result
from Theorem~\ref{thm:abstract}.

We begin with some notation and the setting.
We consider the set $X=\BB^d$ equipped with the measure $d\mu(x)=\omega_\gamma(x)dx$.
Let $Z:= L_0$ and $Y:= L_\cc=L_0+V$ with $V(x)=\cc^2\|x\|^2$,
where we use the notation from above.

%\subsection*{Radial and spherical components of $L_0$ and polynomial basis}

It is not hard to see (e.g. \cite{ZLWZ}) that the operator $L_0$ can be represented as follows
\begin{align}\label{L-0}
L_0 &= - \Delta+\nabla\cdot x(2\gamma-1+x\cdot\nabla)-(2\gamma-1)d \nonumber
\\
& = -\frac{1}{\omega_\gamma(x)}\nabla\cdot \big[\omega_\gamma(x)(1-\|x\|^2)\nabla\big] - \Delta_0
\\
&= -(1-r^2)\partial_r^2-\frac{d-1}{r}\partial_r+(2\gamma+d)r\partial_r -\frac{1}{r^2}\Delta_0, \nonumber
\end{align}
where $\Delta_0$ is the spherical part of the Laplacian $\Delta$.
From above it follows that the operator $L_0$ can be represented as follows (see \cite{KPX2}).
Define
\begin{equation*}
D_{i,j}:=x_i\partial_j-x_j\partial_i,\quad 1 \le i \ne j \le d,
\end{equation*}
and
$$
D^2_{i,i} := [w_\gamma(x)]^{-1} \partial_i \left[(1-\|x\|^2) w_\gamma(x) \right]\partial_i,
\quad 1 \le i \le d.
$$
Then \eqref{L-0} yields (see also (2.2) in \cite{KPX2}):
\begin{equation}\label{decomp}
L_0 = -\sum_{i=1}^d D^2_{i,i} -  \sum_{1\le i < j \le d} D^2_{i,j}
      =  -\sum_{1\le i \le j \le d} D^2_{i,j}.
\end{equation}
Let the domain $D(L_0)$ of the operator $L_0$ be $C^2(\BB^d)$,
which is obviously dense in $L^2(\omega_\gamma)$.
The following basic property of the operator $L_0$ follows from \eqref{decomp} by integration by parts
(see \cite[Theorem~2.1]{KPX2}):
For any $f \in C^2(\BB^d)$ and $g \in C^1(\BB^d)$,
\begin{align}\label{rep-Dmu}
& \int_{\BB^d} L_0 f(x) g(x) w_\mu(x) dx
\\
&= \int_{\BB^d}
\Big[ \sum_{i =1}^d  \partial_i f (x)\partial_i g(x) (1-\|x\|)^2 +
       \sum_{1 \le i< j \le d} D_{i,j} f(x)D_{i,j} g(x) \Big] w_\gamma(x) dx. \notag
\end{align}
Therefore, $L_0$ is a positive symmetric operator.
Moreover, $L_0$ is essentially self-adjoint in $L^2(\omega_\gamma)$.

The associated to $L_0$ quadratic form $\cE$ is defined by
\begin{equation}\label{quad-form-B}
\cE(f, g):= \int_{\BB^d} L_0 f(x) g(x) w_\gamma(x) dx.
\end{equation}
As before we denote its closure by $\overline{\cE}$.
We next show that $\overline{\cE}$ is a Dirichlet form.
Just as in the proof of Theorem~\ref{thm:main} (see also \cite{Fukushima}) for this it suffices to show that
for every $\eps>0$ there exists a function
$\Phi_\eps: \RR\mapsto [-\eps, 1+\eps]$ such that
$\Phi_\eps$ is non-decreasing,
$\Phi_\eps(t)=t$ for $t\in [0, 1]$,
$0\le \Phi_\eps(t') - \Phi_\eps(t) \le t'-t$ if $t<t'$,
and
\begin{equation}\label{BD-cond-2}
f\in D(L_0)
\quad \Longrightarrow\quad
\Phi_\eps(f) \in D(L_0),
\quad
\cE\big(\Phi_\eps(f), \Phi_\eps(f)\big) \le \cE(f, f).
\end{equation}

Let $\Phi_\eps\in C^\infty(\RR)$ have %be an arbitrary function with
the properties: %such that
$-\eps\le \Phi_\eps\le 1+\eps$, $0\le \Phi_\eps'\le 1$ and $\Phi_\eps(t)=t$, $t\in [0, 1]$,
for some $\eps>0$.

Assume $f\in D(L_0)= C^2(\BB^d)$.
Evidently, $\Phi_\eps(f) \in D(L_0)$.
It is easily to see that
$\partial_i \Phi_\eps(f(x))= \Phi_\eps'(f(x))\partial_i f(x)$
and
$D_{ij} \Phi_\eps(f(x))= \Phi_\eps'(f(x)) D_{ij}f(x)$, $i\ne j$.
Hence, using \eqref{rep-Dmu} we get
\begin{align*}
\cE(\Phi_\eps(f),& \Phi_\eps(f))
\\
&=\int_{\BB^d} \Big[ \sum_{i =1}^d  |\partial_i \Phi_\eps(f(x))|^2(1-\|x\|)^2
+ \sum_{1 \le i< j \le d} |D_{i,j} \Phi_\eps(f(x))|^2 \Big] w_\gamma(x) dx
\\
&=\int_{\BB^d} |\Phi_\eps'(f(x))|^2 \Big[ \sum_{i =1}^d  |\partial_i f(x)|^2(1-\|x\|)^2
+ \sum_{1 \le i< j \le d} |D_{i,j}f(x)|^2 \Big] w_\gamma(x) dx
\\
&\le \int_{\BB^d}\Big[ \sum_{i =1}^d  |\partial_i f(x)|^2(1-\|x\|)^2
+ \sum_{1 \le i< j \le d} |D_{i,j}f(x)|^2 \Big] w_\gamma(x) dx
\\
&= \cE(f, f).
\end{align*}
Therefore, $\overline{\cE}$ is a Dirichlet form.
Now, the theory from Section~\ref{sec:abstract} can be applied.
In particular, estimates \eqref{gauss-ball-Lc} follow from Theorem~\ref{thm:abstract} and %\eqref{gauss-ball}.
the Gaussian bounds for $\KK_{\exp(-tL_0)}(x, y)$.
$\qed$

\subsection{H\"{o}lder continuity of the PSWF heat kernel \boldmath $p_t(x, y)$ on $\BB^d$}\label{subsec:Holder-ball}

As is shown in \cite{SS} the quadratic form $\cE$ from \eqref{quad-form-B} generates
a regular strictly local Dirichlet space and the respective Poincar\'{e} inequality is satisfied. % (see \S\ref{subsec:Holder}).
This together with the doubling property of the measure gives another proof of
the Gaussian bounds for $\KK_{\exp(-tL_0)}(x, y)$.
%the Gaussian bounds \eqref{gauss-ball}.
In turn, this %\eqref{gauss-ball}
yields the H\"{o}lder continuity of the kernel $\KK_{\exp(-tL_0)}(x, y)$.
Furthermore, the H\"{o}lder continuity of $\KK_{\exp(-tL_0)}(x, y)$
and \eqref{gauss-ball-Lc}, in light of Theorem~\ref{thm:Holder}, imply
the H\"{o}lder continuous of the PSWF heat kernel $p_t(x,y)$ on $\BB^d$:
There exists a constant $0<\alpha \le 1$ such that
\begin{equation}\label{PSWF-ball-Holder}
|p_t(x, y)-p_t(x', y)|
\le c_1\Big(\frac{\rho(x, x')}{\sqrt{t}}\Big)^\alpha
\frac{\exp\big(- \frac{\rho(x,y)^2}{c_2t}\big)}{\big[V(x, \sqrt t)V(y, \sqrt t)\big]^{1/2}},
\end{equation}
whenever $\rho(x, x')\le \sqrt{t}$.
Here $\alpha$ is the same as for $\KK_{\exp(-tL_0)}(x, y)$.
It is an open question to find the dependence of $\alpha$ on $\gamma$.

\subsection{Besov and Triebel-lizorkin spaces associated to the PSWFs on \boldmath $\BB^d$}\label{subsec:spaces-ball}

The theory of Besov and Triebel-lizorkin spaces on $\BB^d$ associated to the operator $L_\cc$
in the setting described above is now easy to develop
due to the fact that the Gaussian bounds of the heat kernel \eqref{gauss-ball-Lc}
and its H\"{o}lder continuity \eqref{PSWF-ball-Holder} are already established.
These allow to develop the relevant smooth functional calculus and then the associated Besov and Triebel-lizorkin spaces
with a complete set of indices.
To this end one simply follows the general scheme developed in \cite{CKP, KP}.
We will not elaborate on this theory here.


\begin{thebibliography}{99}

\bibitem{ArDe}
    S. Albeverio, Theory of Dirichlet forms and applications,
    In: Lectures on Probability Theory and Statistics, Saint-Flour, 2000,
    Lecture Notes in Math., vol. 1816, pp. 1–106, Springer, Berlin (2003).

 \bibitem{A}
    W. Arendt, M. Demuth, H\"older’s inequality for perturbations of positive semigroups by potentials. J. Math. Anal. Appl. 316 (2006), 652--663.

\bibitem{BCF}
    F. Bernicot, T. Coulhon, and D. Frey,
    Gaussian heat kernel bounds through elliptic Moser iteration,
    J. Math. Pures Appl. (9) 106 (2016), no. 6, 995--1037.

\bibitem{BK}
    A. Bonami, A. Karoui,
    Uniform approximation and explicit estimates for the prolate spheroidal wave functions,
    Constr. Approx.  43 (2016), 15--45.

\bibitem{CKP}
    T. Coulhon, G. Kerkyacharian, P. Petrushev,
    Heat kernel generated frames in the setting of Dirichlet spaces,
    J. Fourier Anal. Appl. 18 (2012), 995--1066.

\bibitem{CS}
    T. Coulhon, A. Sikora,
    Gaussian heat kernel upper bounds via the Phragm\'{e}n-Lindel\"{o}f theorem,
    Proc. London Math. Soc. 96 (2008), 507--544.

\bibitem{DaiXu}
    F. Dai, Y. Xu,
    Approximation theory and harmonic analysis on spheres and balls,
    Springer Monographs in Mathematics, Springer 2013.

\bibitem{Davies-1}
    E. B. Davies,
    Heat kernel and spectral theory,
    Cambridge University Press, 1989.

\bibitem{Davies-2}
    E. B. Davies,
    Linear operators and their spectra,
    Cambridge University Press, 2007.

\bibitem{DKKP}
    S. Dekel, G. Kerkyacharian, G. Kyriazis, and P. Petrushev,
    Compactly supported frames for spaces of distributions associated with nonnegative self-adjoint operators,
    Studia Math. 225 (2014), no. 2, 115--163.

\bibitem{DS}
    N. Dunford, J.T. Schwartz,
    Linear Operators I: General Theory,
    Interscience Publisher, 1958.

\bibitem{DX}
    C. Dunkl, Y. Xu, Orthogonal polynomials of several variables,
    Encyclopedia of Mathematics and its Applications 81, Cambridge University Press, Cambridge, 2001.

\bibitem{DP}
    J. Dziuba\'{n}ski, M. Preisner, Hardy spaces for semigroups with Gaussian bounds,
    Ann. Mat. Pura Appl. (4) 197 (2018), no. 3, 965--987.

\bibitem{Engel-Nagel}
    K.-J. Engel, R. Nagel,
    A short course on operator semigroups, Universitext, Springer, New York, 2006.

\bibitem{Fukushima}
    M. Fukushima, Y. Oshima, M. Takeda,
    Dirichlet forms and symmetric Markov processes,
    %de Gruyter Studies in Math., vol. 19, Walter de Gruyter, Berlin and Hawthorne, NY, 1994.
    De Gruyter Studies in Mathematics, 19, Walter de Gruyter \& Co., Berlin, 1994.

\bibitem{F}
    C. Flammer,
    Spheroidal Wave Functions, Stanford University Press, Stanford, CA, 1957.

\bibitem{FJ1}
    M. Frazier, B. Jawerth,
    Decomposition of  Besov Spaces,
    Indiana Univ. Math. J. 34 (1985), 777--799.

\bibitem{FJ2}
    M. Frazier, B. Jawerth,
    A discrete transform and decompositions of distribution spaces,
    J. Funct. Anal. 93 (1990), 34--170.

\bibitem{FJW}
    M. Frazier, B. Jawerth, and G. Weiss,
    Littlewood-Paley Theory and the Study of Function  Spaces,
    CBMS 79 (1991), AMS.

\bibitem{GS}
    P. Gyrya and L. Saloff-Coste,
    Neumann and Dirichlet heat kernels in inner uniform domains, Ast\'{e}risque (2011), no. 336, viii+144.

\bibitem{HS}
    W. Hebisch, L. Saloff-Coste, On the relation between elliptic and parabolic Harnack inequalities,
    Ann. Inst. Fourier (Grenoble) 51 (2001), no. 5, 1437--1481.

\bibitem{HL}
    J.A. Hogan, J.D. Lakey,
    Duration and Bandwidth Limiting, Prolate Functions, Sampling, and Applications,
    Applied and Numerical Harmonic Analysis, Birkh\"{a}user/Springer, New York, 2012.

\bibitem{KP}
    G. Kerkyacharian, P. Petrushev,
    Heat kernel based decomposition of spaces of distributions in the framework of Dirichlet spaces.
    Trans. Amer. Math. Soc. 367 (2015), 121--189.

\bibitem{KPX1}
    G. Kerkyacharian, P. Petrushev, Y. Xu,
    Gaussian bounds for the weighted heat kernels on the interval, ball, and simplex,
    Constr. Approx. 51 (2020), no. 1, 73--122.

\bibitem{KPX2}
    G. Kerkyacharian, P. Petrushev, Y. Xu,
    Gaussian bounds for the heat kernels on the ball and the simplex: classical approach,
    Studia Math. 250 (2020), no. 3, 235--252.

\bibitem{MF}
    P. M. Morse, H. Feshbach,
    Methods of Theoretical Physics, Two volumes, McGraw-Hill Book Co., Inc., New York, 1953.

\bibitem{ORX}
    A. Osipov, V. Rokhlin, H. Xiao,
    Prolate Spheroidal Wave Functions of Order Zero.
    Mathematical tools for bandlimited approximation, Applied Mathematical Sciences, 187, Springer, New York, 2013.

    \bibitem{Ouhabaz}
    E. M. Ouhabaz, Analysis of Heat Equations on Domains, Princeton Univ. Press, Princeton, NJ, 2005.

\bibitem{PX}
    P. Petruushev, Y. Xu,
    Localized polynomial frames on the interval with Jacobi weights,
    J. Fourier Anal. and Appl. 11 (2005), 557--575.

\bibitem{SC}
    L. Saloff-Coste, Aspects of Sobolev-type inequalities,
    London Mathematical Society Lecture Note Series, vol. 289, Cambridge University Press, Cambridge, 2002.

\bibitem{SlepP}
    D. Slepian, H. O. Pollak,
    Prolate spheroidal wave functions, Fourier analysis and uncertainty I,
    Bell System Tech. J. 40 (1961), 43--64.

\bibitem{Slep-1}
    D. Slepian,
    Prolate spheroidal wave functions,
    Fourier analysis and uncertainty--IV: Extensions to many dimensions; generalized prolate spheroidal functions,
    Bell System Tech. J. 43 (1964), 3009--3057.

\bibitem{Slep-2}
    D. Slepian,
    Some Asymptotic Expansions for Prolate Spheroidal Wave Functions,
    J. Math. Phys. 44, No. 2, (1965), 99--140.

\bibitem{Sz}
    G. Szeg\"o, Orthogonal polynomials,
    Amer. Math. Soc. Colloq. Publ. Vol. 23, Amer. Math. Soc.
    Providence, 1975.

\bibitem{SS}
    P. Sj\"{o}gren, T. Szarek,
    Analysis in the multi-dimensional ball,
    Mathematika 65 (2019), no. 2, 190--212.

\bibitem{Wang_Zhang}
    L.-L. Wang and J. Zhang,
    A new generalization of the PSWFs with applications to spectral approximations on quasi-uniform grids,
    Appl. Comput. Harmon. Anal. 29, (2010), 303--329.

\bibitem{ZLWZ}
    J. Zhang, H. Li, Li-Lian Wang, Z. Zhang,
    Ball prolate spheroidal wave functions in arbitrary dimensions,
    Appl. Comput. Harmon. Anal. 48 (2020), 539--569.
\end{thebibliography}
\end{document}